\documentclass[11pt]{article} 

\usepackage[margin=1in]{geometry}
\usepackage{setspace}

\normalsize
\usepackage{amsfonts}
\usepackage{amssymb}
\usepackage[cmex10]{mathtools}
\usepackage{amscd}
\usepackage{graphics}
\usepackage[dvips]{graphicx}
\usepackage{epsfig}
\usepackage{array}
\usepackage{eqparbox}
\usepackage{hyperref}
\usepackage{fancyhdr}
\usepackage{appendix}
\usepackage{titling}
\usepackage{bm}
\usepackage{enumerate}
\usepackage{tabularx}
\usepackage[super]{nth}
\usepackage{caption} 
\usepackage{subcaption}

\usepackage{booktabs}
\usepackage{float} 

\usepackage[round]{natbib}  
\DeclareMathAlphabet{\mathbsf}{OT1}{cmss}{bx}{n}
\DeclareMathAlphabet{\mathssf}{OT1}{cmss}{m}{sl}
\DeclareMathAlphabet{\mathcsf}{OT1}{cmss}{sbc}{n}
\DeclareMathOperator{\arccot}{arccot}

\usepackage{stackengine}
\newcommand{\ie}{{\em i.e.}}

\newcommand{\eg}{{\em e.g.}}

\newcommand{\secref}[1]{Section~\ref{#1}}
\newcommand{\figref}[1]{Fig.~\ref{#1}}
\newcommand{\tabref}[1]{Table~\ref{#1}}

\newcommand{\keywords}[1]{\textbf{Keywords:} #1}

\makeatletter
\def\blfootnote{\xdef\@thefnmark{}\@footnotetext}
\makeatother

\newcommand{\qed}{\nobreak \ifvmode \relax \else
      \ifdim\lastskip<1.5em \hskip-\lastskip
      \hskip1.5em plus0em minus0.5em \fi \nobreak
      \vrule height0.75em width0.5em depth0.25em\fi}

\usepackage{newlfont}
\date{}
\begin{document}
\title{\Large{\textbf{A Laplace equation approach to the Behrens--Fisher problem}}}
\author{Nagananda K G and Jong Sung Kim\thanks{The authors are with Fariborz Maseeh Department of Mathematics and Statistics, Portland State University, Portland, OR 97201, USA. E-mail: \texttt{\{nanda, jong\}@pdx.edu}.}}
\setlength{\droptitle}{-1.in}
\maketitle
\vspace{-2cm}

\begin{abstract}
We develop a partial differential equation formulation of the Behrens--Fisher problem for two independent normal samples with unknown and unequal variances. An orthogonal decomposition separates mean and residual components (corresponding to the centered within-sample variation left after removal of the mean directions) and recasts the studentized difference of sample means as a scale-invariant geometric constraint. This reduction transforms the distributional problem into the evaluation of spherical wedge probabilities, which are identified with harmonic measure and with the value at the origin of a Laplace--Dirichlet boundary value problem. From this framework, we derive exact finite-sample representations for the cumulative distribution function and the probability density function in terms of beta functions, with dependence only on the sample sizes and the variance ratio. These representations place the Behrens--Fisher law in a standard special-function form that is directly accessible in widely available commercial software\textemdash including Microsoft Excel\textemdash thereby facilitating distributional evaluation and quantile computation. We also obtain a Gegenbauer separation-of-variables expansion for the associated harmonic extension and its threshold derivative, with coefficients in closed Beta--Gamma form, and derive sharp tail expansions with explicit leading constants and higher-order corrections.
\end{abstract}

\section{Introduction}\label{sec:introduction}\vspace{-0.3cm}
Consider two independent random samples $\left(X_1, \dots, X_{n_1}\right)$ and $\left(Y_1, \dots, Y_{n_2}\right) $, drawn from two normal populations with means $\mu_1, \mu_2$ and variances $\sigma_1^2, \sigma_2^2$, respectively. For the two samples, let $\nu_i = n_i - 1$, $i=1, 2$, denote the usual degrees of freedom associated with the sample variances.  The corresponding sample means are $\bar{X} = \sum_{i=1}^{n_1} X_i/n_1$ and $\bar{Y} = \sum_{i=1}^{n_2} Y_i/n_2$, and the sample variances are $S_1^2 = \sum_{i=1}^{n_1} (X_i - \bar{X})^2/\nu_1$ and  $S_2^2 = \sum_{i=1}^{n_2} (Y_i - \bar{Y})^2/\nu_2$.  To simplify later expressions, we define $g = \sigma_1^2/n_1$, $h = \sigma_2^2/n_2$, $\kappa = \sigma_1^2/\sigma_2^2$, where $g$ and $h$ are the population variance contributions to the sampling variances of $\bar{X}$ and $\bar{Y}$, and $\kappa$ is the ratio of the two population variances.  We are interested in testing the null hypothesis $H_0 : \mu_1 = \mu_2$, when both population variances $(\sigma_1^2, \sigma_2^2)$ are unknown and are not assumed equal, that is, when $\sigma_1^2 \neq \sigma_2^2$.  This is the classical Behrens--Fisher problem, dating back to the work of \cite{Behrens1929} and \cite{Fisher1935}.  We first present the central mathematical challenge posed by this problem, and then outline our approach in \secref{subsec:main_contribution}.

The Behrens--Fisher test statistic is given by 
\begin{eqnarray}
T = \frac{\bar{X} - \bar{Y}}{\sqrt{\frac{S_1^2}{n_1} + \frac{S_2^2}{n_2}}}. 
\label{eq:Behrens_Fisher_statistic}
\end{eqnarray}
Under the null hypothesis $H_0$, the difference of the sample means satisfies $\bar{X} - \bar{Y} \sim \mathcal{N}(0, g + h)$.  In addition, the sample variances satisfy $\nu_1 S_1^2/\sigma_1^2 \sim \chi_{\nu_1}^2$ and $\nu_2 S_2^2/\sigma_2^2 \sim \chi_{\nu_2}^2$, and these chi-square variables are independent of $\bar{X}$ and $\bar{Y}$.  The main difficulty arises from the denominator of $T$. Note that, the denominator in \eqref{eq:Behrens_Fisher_statistic} is not a simple pooled estimator, but rather a linear combination of two independent scaled chi-square random variables involving the unknown variances $(\sigma_1^2, \sigma_2^2)$.  Furthermore, when the population variances are unequal ($\sigma_1^2 \neq \sigma_2^2$), this structure prevents the statistic from reducing to the usual Student $t$ form (see \cite[Exercise 8.42]{Casella2001}). Define $W_i = \nu_i S_i^2/\sigma_i^2$,  $ i=1, 2$.  Then, the statistic $T$ can be rewritten as
\begin{eqnarray*}
T &=& \frac{\bar{X} - \bar{Y}}{\sqrt{\frac{g W_1}{\nu_1} + \frac{h W_2}{\nu_2}}}.
\end{eqnarray*}
Hence, conditioned on $\sigma_1^2$ and $\sigma_2^2$, the distribution of $T$ is given by 
\begin{eqnarray*}
T \mid \sigma_1^2, \sigma_2^2 &\sim& \frac{\mathcal{N}(0, g+h)}{\sqrt{\frac{g W_1}{\nu_1} + \frac{h W_2}{\nu_2}}}.
\end{eqnarray*}
This should be contrasted with the equal-variance case. If $\sigma_1^2 = \sigma_2^2 = \sigma^2$ and the common variance $\sigma^2$ is unknown, then one may use the pooled variance estimator
\begin{eqnarray}
S_p^2 &=& \frac{1}{\nu_1 + \nu_2} \left( \sum_{j=1}^{n_1} (X_j - \bar{X})^2 + \sum_{j=1}^{n_2} (Y_j - \bar{Y})^2 \right).
\label{eq:pooled_estimator}
\end{eqnarray}
In that case, the corresponding likelihood-ratio-type statistic \eqref{eq:Behrens_Fisher_statistic} becomes
\begin{eqnarray*}
T \mid \sigma^2 &=& \frac{\bar{X} - \bar{Y}}{\sqrt{S_p^2 \left( \frac{1}{n_1} + \frac{1}{n_2} \right)}},
\end{eqnarray*}
and this statistic has the Student $t$ distribution with $\nu_1 + \nu_2$ degrees of freedom (see \cite[Exercise 8.41]{Casella2001}): $T \mid \sigma^2 \sim t_{\nu_1 + \nu_2}$.  Thus, when the population variances are equal, the classical two-sample $t$-test applies directly. 

By contrast, when $\sigma_1^2 \neq \sigma_2^2$ and $(\sigma_1^2, \sigma_2^2)$ is unknown, no analogous simplification is available. In particular, the pooled estimator $S_p^2$ in \eqref{eq:pooled_estimator} is no longer appropriate, and the distribution of
\begin{eqnarray*}
T &=& \frac{\bar{X} - \bar{Y}}{\sqrt{\frac{g W_1}{\nu_1} + \frac{h W_2}{\nu_2}}}
\end{eqnarray*}
must instead be obtained by integrating over the joint distribution of $W_1$ and $W_2$. Accordingly, 
\begin{eqnarray*}
f_T(t) &=& \int_0^\infty \int_0^\infty f_{T \mid W_1, W_2}(t) \, f_{W_1}(w_1) \, f_{W_2}(w_2) \, dw_1 \, dw_2,
\end{eqnarray*} 
where $f_{W_1}$ and $f_{W_2}$ are the chi-square densities corresponding to $W_1$ and $W_2$.  Conditioned on $W_1 = w_1$ and $W_2 = w_2$, the random variable $T$ is normally distributed with mean zero and variance determined by $w_1$ and $w_2$. Therefore, its marginal density can be written as
\begin{eqnarray}
f_T(t) = \int_0^\infty \int_0^\infty \frac{\sqrt{g w_1 / \nu_1 + h w_2 / \nu_2}}{\sqrt{2\pi (g+h)}} \exp \left( -\frac{t^2 (g w_1 / \nu_1 + h w_2 / \nu_2)}{2(g+h)} \right) f_{W_1}(w_1) f_{W_2}(w_2) \, dw_1 \, dw_2.
\label{eq:challenge}
\end{eqnarray}

Over several decades, this problem has attracted sustained attention, leading to a broad literature on approximations, modifications, and practically implementable procedures; see, for example, the results of  \cite{Welch1947,Scheffe1970,Bozdogan1986,Best1987,Asiribo1989,Duong1996,Kim1998a,Vangel2005,Dudewicz2007,Chang2008,Nadarajah2017,Chaturvedi2019,Wang2022,Chen2022,Chen2023}. Related developments have also appeared in nonparametric formulations; see, for example, the work of \cite{Brunner2002, Larocque2010, Konietschke2012, Konietschke2012a, He2022}.  Considerable work has also focussed on high-dimensional extensions; see, for instance, \cite{Zhou2017,Cao2019,Zhang2021,Pei2026}.  Recent investigation by \cite{Nagananda2026} collapsed the challenging two-dimensional integral in \eqref{eq:challenge} to a tractable single-contour integral.  

Despite consistent efforts, evaluating the integral in \eqref{eq:challenge} does not lead to a closed-form expression for the cumulative distribution function (cdf) and, hence, the probability density function (pdf) of the statistic in \eqref{eq:Behrens_Fisher_statistic} in terms of standard named distributions.  This is the essential analytical difficulty of the Behrens--Fisher problem: once the equal-variance assumption is dropped, the familiar two-sample $t$ framework no longer yields an exact distributional formula of the usual elementary kind.  That lack of an explicit representation makes exact significance calculations, power analysis, and precise error assessment substantially more difficult, especially in small samples or when the two variances differ markedly.  We aim to resolve this issue by taking a completely new standpoint. 

\vspace{-0.3cm}
\subsection{Main contribution}\label{subsec:main_contribution} \vspace{-0.3cm}
We employ a partial differential equation (PDE)\textemdash specifically, the Laplace's equation (see \cite[Chapter 4]{McOwen2003})\textemdash approach to derive the sampling distribution of the statistic \eqref{eq:Behrens_Fisher_statistic}.  This approach transforms a ratio of random quantities with unknown, unequal variances into a clean, geometric boundary-value problem.  After the right orthogonal change of variables, the studentized inequality $T \leq t$ is a scale-invariant statement about a direction on the unit sphere, not about radial magnitude. That makes $\mathbb{P}(T \leq t)$ a question about how much of the sphere is cut out by a spherical wedge. In this setting, the value at the origin of the harmonic (Laplacian) solution with $0/1$ boundary data equals that probability (harmonic measure). This approach eliminates ad hoc approximations and tedious conditioning arguments; the distribution emerges from a canonical PDE construct.

The method proceeds by an orthogonal decomposition that isolates mean and residual components, the latter corresponding to the centered within-sample variation left after removal of the mean directions. This is followed by a geometric re-expression of the event $\{T \le t\}$ as a scale-invariant conic constraint. The distributional problem is thereby reduced to computing the fraction of directions on a high-dimensional sphere that lie in a corresponding spherical wedge, {\ie}, the normalized spherical surface measure of the set of directions for which the event defining the cdf occurs. We then show that this probability is equal to the value at the origin of a harmonic function solving a Laplace--Dirichlet boundary value problem \cite[Chapter 4.2e]{McOwen2003}, thereby linking the cdf directly to harmonic measure.  Exploiting this structure yields canonical expressions for the cdf and the pdf given by \eqref{eq:main_result_cdf} and \eqref{eq:main_result_pdf}, respectively: 
\begin{eqnarray}
\text{cdf}:~F_T(t) &=& I_{\frac{1}{\left(1+t^2\right)}}\left(\frac{D-1}{2}, \frac{1}{2}\right), \label{eq:main_result_cdf} \\
\text{pdf}:~f_{T}(t) &=& \frac{1}{c_{\kappa} B\left(\frac{1}{2}, \frac{D-1}{2}\right)} \left(1+\frac{t^{2}}{c_{\kappa}^{2}}\right)^{-\frac{D}{2}}, \label{eq:main_result_pdf}
\end{eqnarray}
where $c_{\kappa} = \sqrt{\kappa/n_{1} + 1/n_{2}}$, the variance ratio $\kappa = \sigma_1^2/\sigma_2^2$, $D = n_{1}+n_{2}-1$, $B(a, b)$ is the beta function explained systematically by \cite{Paris2010}, and $I_x(a, b)$ is the regularized incomplete beta function treated in depth by \cite{Zelen1972}.  Both \eqref{eq:main_result_cdf} and \eqref{eq:main_result_pdf} are readily evaluable using commercially available software\textemdash for instance, Microsoft Excel\textemdash and depend only on the sample sizes $(n_1, n_2)$ and the variance ratio $\kappa$, but not on the individual variances. 

The derivations of \eqref{eq:main_result_cdf} and \eqref{eq:main_result_pdf} proceed via the following steps:
\begin{enumerate}[(a)]
\item \textbf{Orthogonal decomposition:} We being by splitting each sample into two non-overlapping pieces: a single mean direction (that drives the difference of sample means) and a residual cloud (that drives the sample variances). This rotation makes the mean part independent of the residual part, so the numerator and denominator of Behrens--Fisher's $T$ in \eqref{eq:Behrens_Fisher_statistic} come from independent sources.  This is detailed in \secref{sec:orthogonal_decomposition}.  After separating the mean and residual coordinates, an appropriate rescaling converts the event $\{T\le t\}$ into a scale-invariant conic inequality. Thus the problem depends only on direction\textemdash reducing the distributional question to spherical geometry and hence to a Laplace--Dirichlet boundary value problem.

\item \textbf{Cone (positive-homogeneous) inequality:} After the orthogonal decomposition and rescaling of the event $\{T\le t\}$, the statistic $T$ is expressed in the new coordinates $(u, v, w)$, where $u$ is the scalar coordinate in the mean-difference direction and $v$ and $w$ are the residual coordinate vectors for the two samples. Consequently, the event $\{T \le t\}$ is identified with the subset of the transformed sample space consisting of all $(u,v,w)$ satisfying $u \le t\sqrt{\lVert v \rVert^2 + \lVert w \rVert^2}$.  Since the defining inequality is positively homogeneous, this subset is invariant under multiplication by positive scalars and therefore forms a dihedral cone in the coordinates $(u, v, w)$.  Intersecting this cone with the unit sphere produces a spherical wedge $W_t$, that is, the subset of the sphere consisting of those unit vectors that satisfy the same inequality, and hence represent the admissible directions corresponding to the event $\{T \le t\}$.  After the orthogonal rotation and rescaling, the event $\{T \le t\}$ is equivalent to the statement that the mean coordinate is no larger than $t$ times the Euclidean size of the residual components. Since the defining inequality is positively homogeneous, $\{T \le t\}$ is invariant under radial scaling, and its probability depends only on direction. The problem thereby reduces to computing the normalized spherical surface measure of $W_t$.  This development can be found in \secref{sec:cone_inequality}.

\item \textbf{PDE (Laplace's equation) that returns the cdf:}  In \secref{sec:Laplace_equation_cdf}, we solve the Laplace's equation $\Delta u=0$ in the unit ball with boundary data $1$ on $W_t$ and $0$ elsewhere; then $F_T(t)=u(0)$ by harmonic measure. For the spherical wedge arising here, this harmonic measure admits a closed-form evaluation in terms of a regularized incomplete beta function, so the PDE formulation leads directly to an explicit expression for the cdf. On the unit ball, we therefore assign boundary value $1$ on the wedge and $0$ elsewhere, and solve Laplace's equation in the interior. The value of the resulting harmonic solution at the center is exactly the cdf $F_T(t)=\mathbb P(T\le t)$. Thus, the sampling-distribution problem for the statistic in \eqref{eq:Behrens_Fisher_statistic} is reduced to evaluating a canonical boundary-value solution at one point, bringing to bear standard analytic tools such as harmonic measure, the Poisson kernel, and eigenfunction expansions.  

\item \textbf{Separation of variables on the spherical wedge:} In \secref{sec:separation_of_variables}, we express the Laplacian in spherical coordinates adapted to the wedge and decompose the boundary value problem into radial and angular parts. The angular equation reduces to a classical Sturm--Liouville problem on an interval (see \cite[Chapter 4,  pp. 105]{McOwen2003}), with boundary conditions inherited from the wedge geometry, and its eigenfunctions are given by standard Gegenbauer families treated by \cite{Durand1976}.  The  full analysis can be found in \secref{sec:Gegenbauer_coefficients}. This yields a convergent eigenfunction expansion for the harmonic solution and, after evaluation at the origin and differentiation with respect to a threshold, provides a series representation for the pdf (see \secref{subsec:expression_pdf}).  Beyond exact distributional formulas, the Gegenbauer expansion also identifies the leading angular mode that governs the tail behaviour, thereby producing precise large-threshold tail expansions with explicit sharp tail constants and a systematic hierarchy of higher-order correction terms.  Relevant details are in \secref{subsec:tail_behavior}.  Some graphical results and numerical tables are presented in \secref{subsec:numerical_results} to complement the theoretical development. 
\end{enumerate}
The novelty of our contribution is twofold.  First, the approach gives a new analytic derivation of the exact finite-sample distribution in the Behrens--Fisher setting with unknown and unequal variances. Rather than proceeding through classical integral representations or characteristic-function arguments, the distribution is obtained from a Laplace--Dirichlet problem on a spherical wedge via separation of variables. This results in a convergent Gegenbauer eigenfunction expansion, while the cdf and pdf emerge naturally in terms of well-known beta functions.  The same framework also leads to explicit tail asymptotics with sharp leading constants.  Second, the geometric formulation makes extension natural: replacing the Euclidean Laplacian by a suitable anisotropic divergence-form operator suggests corresponding analogues for elliptical families, while preserving the same boundary-value perspective.

One of the most widely used methods in practice to address the Behrens--Fisher problem was derived by \cite{Welch1947}.   Welch's test yields a practically excellent approximation by comparing the statistic to a Student's $t$ distribution with an estimated effective degrees of freedom
\begin{eqnarray*}
\hat{\nu} \approx  \frac{\left(\frac{S_1^2}{n_1} + \frac{S_2^2}{n_2}\right)^2}{\left[\frac{S_1^4}{n^2_1\nu_1} + \frac{S_2^4}{n^2_2\nu_2}\right]}, 
\end{eqnarray*}
obtained using the approximation developed by \cite{Satterthwaite1946}. Its strength lies in its robustness and good Type I error control without requiring $\sigma_1^2 = \sigma_2^2$ (see \cite[Chapter 8]{Casella2001}).  Our result, however, is \emph{not} another approximation of that type.  Instead, we provide an exact finite-sample characterization of the distribution of the Behrens--Fisher statistic by recasting the event $\{T \le t\}$ as a geometric constraint and identifying the resulting probability with the value at the origin of a Laplace--Dirichlet problem.  This difference brings three key advantages: (i) The expressions derived in this paper isolate the precise dependence of the law on the sample sizes and the variance ratio, rather than compressing that structure into a single approximate degrees-of-freedom parameter. (ii) The representation in terms of beta functions gives a standard special-function form that is both interpretable and numerically accessible in common software.  (iii) The PDE and harmonic-measure framework supplies additional analytic information that Welch's approximation does not provide naturally, including a compact density representation, a separation-of-variables expansion, and sharp tail asymptotics with explicit constants.  

\section{Orthogonal decomposition}\label{sec:orthogonal_decomposition} 
Orthogonal decomposition is a systematic way to separate what matters for the Behrens--Fisher statistic from what does not\textemdash and to make independence visible. Starting with the raw samples, we project each vector onto the 1-dimensional mean direction and its $(n_i-1)$-dimensional residual subspace. This yields two coordinates per group: a mean component (driving $\bar X-\bar Y$) and an orthogonal residual block (driving $S_1^2$ or $S_2^2$). By rotational invariance of the multivariate normal and Cochran's theorem (see \cite{Cochran1934}), the mean coordinates are independent of the residual coordinates, and the residual sums of squares become simple chi-square norms. In the unknown--unequal variance regime, the decomposition separates the mean-contrast (numerator) and the residual-based studentizer (denominator) into orthogonal\textemdash and therefore independent\textemdash subspaces, so the statistic $T$ in \eqref{eq:Behrens_Fisher_statistic} becomes a ratio of functions of independent components.

The orthogonal decomposition is the pivotal step that leads to a geometric\textemdash and ultimately PDE-based\textemdash treatment.  After standardizing by the unknown scales through the variance ratio $\kappa$, the studentized inequality $T\le t$ becomes a positive-homogeneous relation comparing a single mean-contrast coordinate to the Euclidean norm of the residual blocks.  In the orthogonal coordinates, this is a cone $u \le t \sqrt{\lVert v \rVert^2 + \lVert w \rVert^2}$, so the event depends only on direction, not radius.  This observation turns the probability $\mathbb P(T \le t)$ into a question about the portion of the unit sphere cut out by a spherical wedge\textemdash precisely the setting where harmonic measure and Laplace's equation apply.

Once in this form, the problem becomes the Dirichlet problem for Laplace's equation on a cone-capped ball: the cdf equals the harmonic measure of the spherical cap, and separation of variables on the spherical wedge gives an explicit eigenfunction expansion in Gegenbauer polynomials.  This leads to closed-form Beta--Gamma formulas for Gegenbauer coefficients, exact expressions for the cdf/pdf, and sharp tail asymptotics governed by the principal spherical eigenvalue.  The same framework also opens avenues beyond the classical normal case ({\eg}, elliptical laws via anisotropic Laplacians) and to other statistics whose rejection regions become conic after the same orthogonalization.

In the following, we provide full details of expressing $T$ in orthogonal coordinates that separate mean and residual components and make independence/chi-square structure explicit.  Fix integers $n_1\ge 2$, $n_2\ge 2$.  Let $1_m\in\mathbb R^m$ denote the $m \times 1$ vector of all ones.  For a matrix $A$, $A^\top$ denotes the transpose; $I_m$ is the $m\times m$ identity.  For vectors $x$, $\lVert x \rVert = \sqrt{x^\top x}$ is the Euclidean norm.  The observations are denoted by $X = (X_{1}, \dots, X_{n_{1}})^\top$ and $Y = (Y_{1}, \dots, Y_{n_{2}})^\top$.  The distributional assumptions under the null $H_0:\mu_1=\mu_2=\mu$ are $X \sim \mathcal N(\mu1_{n_1}, \sigma_1^2 I_{n_1})$, $Y \sim \mathcal N(\mu1_{n_2}, \sigma_2^2 I_{n_2})$,  with unknown $(\sigma_1^2, \sigma_2^2)$ and $\sigma_1^2\neq \sigma_2^2$. Also, $X$ and $Y$ are independent.  

\subsection{Projection operators and orthogonal bases}
We first separate the mean direction $\mathrm{span}{\{1_m\}}$ from its orthogonal complement in each sample.  The rank-1 mean projector in $\mathbb R^m$ is defined by $P_m = \frac{1}{m}1_m 1_m^\top$.   Then, $P_m x = \left(\frac{1}{m}1_m^\top x\right)1_m$.  The residual (centering) projector is defined by $M_m = I_m - P_m$.  Then, $M_m x = x - (\text{average of entries of }x)\times 1_m$, so $M_m x$ has components summing to zero.  The following basic identities follow directly:  $P_m^2=P_m$,  $M_m^2 = M_m$, $P_m M_m = 0$, $P_m + M_m = I_m$, $\operatorname{rank}(P_m)=1$, $\operatorname{rank}(M_m) = m-1$.  Any orthogonal matrix $U_1\in\mathbb R^{n_1 \times n_1}$ whose first column is $u_1 = \frac{1_{n_1}}{\sqrt{n_1}}$ can be chosen, and $u_1$ can be completed to an orthonormal basis $\{u_1,u_2,\dots,u_{n_1}\}$ of $\mathbb R^{n_1}$ such that $\operatorname{span}\{u_2,\dots,u_{n_1}\} = \operatorname{Im}(M_{n_1})={\text{vectors orthogonal to }1_{n_1}}$.  Thus, $U_1^\top U_1 = I_{n_1}$ and $\{u_2,\dots,u_{n_1}\}$ is an orthonormal basis for the residual subspace.  Similarly, choose $U_2 \in \mathbb R^{n_2\times n_2}$ orthogonal with $v_1 = \frac{1_{n_2}}{\sqrt{n_2}}$ as the first column, and $\{v_2, \dots, v_{n_2}\}$ an orthonormal basis for $\operatorname{Im}(M_{n_2})$.

We rotate $X-\mu1_{n_1}$ and $Y-\mu1_{n_2}$ into mean and residual components.  The rotated coordinates for $X$ and $Y$ are defined by 
\begin{eqnarray}
\begin{pmatrix} m_X\\  r_X \end{pmatrix}
&=& U_1^\top \left(X - \mu1_{n_1}\right), \quad
m_X\in\mathbb R,\ \ r_X\in\mathbb R^{n_1-1}, \label{eq:rotated_coordinates_X} \\
\begin{pmatrix} m_Y \\ r_Y \end{pmatrix}
&=& U_2^\top \left(Y-\mu1_{n_2}\right),\quad
m_Y\in\mathbb R,\ \ r_Y\in\mathbb R^{n_2-1}, \label{eq:rotated_coordinates_Y}
\end{eqnarray}
where $m_X$ and $m_Y$ are given by  
\begin{eqnarray*}
m_X &=& u_1^\top (X - \mu1_{n_1}) = \frac{1}{\sqrt{n_1}} \sum_{i=1}^{n_1} (X_i - \mu) = \sqrt{n_1}(\bar{X} - \mu), \\
m_Y &=&  v_1^\top (Y - \mu1_{n_2}) = \frac{1}{\sqrt{n_2}} \sum_{i=1}^{n_2} (Y_i - \mu) = \sqrt{n_2}(\bar{Y} - \mu),
\end{eqnarray*}
since the first row of $U_1^\top$ is $u_1^\top = \frac{1}{\sqrt{n_1}}1_{n_1}^\top$.  The residual blocks are given by 
\begin{eqnarray*}
r_X = \begin{pmatrix}u_2^\top\\\ \vdots\\ u_{n_1}^\top\end{pmatrix}(X - \mu1_{n_1}), \qquad r_Y = \begin{pmatrix}v_2^\top\\ \vdots\\ v_{n_2}^\top\end{pmatrix} (Y - \mu1_{n_2}).
\end{eqnarray*}
Since $U_1$ is orthogonal, $\lVert X-\mu1_{n_1}\rVert^2 = m_X^2 + \lVert r_X \rVert^2$. Furthermore, since $M_{n_1}$ is idempotent and $M_{n_1} = I_{n_1}-\frac{1}{n_1}1_{n_1}1_{n_1}^\top$ , we have
\begin{eqnarray*}
\lVert r_X \rVert^2 &=& \left\lVert (U_1^\top (X-\mu1_{n_1}))_{2:n_1}\right\rVert^2 = (X-\mu1{n_1})^\top M_{n_1}(X-\mu1_{n_1}) \nonumber \\
&=& \sum_{i=1}^{n_1}(X_i-\mu)^2 - \frac{1}{n_1}\left(\sum_{i=1}^{n_1} (X_i-\mu)\right)^2 = \sum_{i=1}^{n_1}(X_i - \bar{X})^2 \nonumber \\
&=& (n_1-1)S_1^2.
\end{eqnarray*}
An identical argument yields $\lVert r_Y \rVert^2 = \sum_{j=1}^{n_2}(Y_j-\bar{Y})^2=(n_2-1)S_2^2$.

We can now express the distributions of the rotated coordinates.  Since $X\sim \mathcal \mathcal{N}(\mu1_{n_1},\sigma_1^2 I_{n_1})$ and $U_1$ is fixed orthogonal, rotational invariance of the multivariate normal yields $U_1^\top (X - \mu1_{n_1}) \sim \mathcal N\left(0,\ \sigma_1^2 I_{n_1}\right)$. Therefore, we have $m_X \sim \mathcal N(0,\sigma_1^2)$, and $r_X \sim \mathcal N\left(0,\ \sigma_1^2 I_{n_1 - 1}\right)$.  Moreover, since the covariance is diagonal $\sigma_1^2 I_{n_1}$ in the rotated basis, the scalar coordinate $m_X$ and the vector $r_X$ are independent.  Similarly, we can write $m_Y \sim \mathcal N(0,\sigma_2^2)$, and $r_Y \sim \mathcal N\left(0,\ \sigma_2^2 I_{n_2 - 1}\right)$, and $m_Y$ and $r_Y$ are independent. Also, since $X \perp Y$, the four blocks $(m_X,r_X)$ and $(m_Y,r_Y)$ are mutually independent.  This results in the following chi-squared structure, independently of $m_X,m_Y$:
\begin{eqnarray}
\frac{\lVert r_X \rVert^2}{\sigma_1^2} &=& \sum_{k=2}^{n_1} \left(\frac{u_k^\top (X - \mu1_{n_1})}{\sigma_1}\right)^2 \sim \chi^2_{n_1 - 1}, \label{eq:chi_rX} \\
\frac{\lVert r_Y \rVert^2}{\sigma_2^2} &=& \sum_{k=2}^{n_2} \left(\frac{v_k^\top (Y - \mu1_{n_2})}{\sigma_2}\right)^2 \sim \chi^2_{n_2 - 1}. \label{eq:chi_rY}
\end{eqnarray}

The aforementioned development lets us write the statistic $T$ in the desired rotated coordinates.  Firstly, the sample means and sample variances are given by
\begin{eqnarray*}
\bar X = \mu+\frac{m_X}{\sqrt{n_1}}, &\quad& \bar Y=\mu+\frac{m_Y}{\sqrt{n_2}}, \\
\sum_{i=1}^{n_1}(X_i-\bar X)^2 = \lVert r_X \rVert^2, &\quad& S_1^2 = \frac{\lVert r_X \rVert^2}{n_1-1}, \\
\sum_{j=1}^{n_2}(Y_j-\bar Y)^2 = \lVert r_Y \rVert^2, &\quad& S_2^2 = \frac{\lVert r_Y \rVert^2}{n_2-1}.
\end{eqnarray*}
To expose spherical symmetry, we standardize the mean and residual coordinates by the true scales.  Towards this end, the standardized mean and residual coordinates are defined as follows: 
\begin{eqnarray*}
\tilde m_X &=& \frac{m_X}{\sigma_1} \in \mathbb R, \quad \tilde{r}_X = \frac{r_X}{\sigma_1}\in\mathbb R^{n_1-1}, \\
\tilde m_Y &=& \frac{m_Y}{\sigma_2} \in \mathbb R,\quad \tilde{r}_Y = \frac{r_Y}{\sigma_2}\in\mathbb R^{n_2-1},
\end{eqnarray*}
whose distributions are given by $\tilde m_X \sim \mathcal N(0,1)$, $\tilde m_Y \sim \mathcal N(0,1)$, $\tilde{r}_X\sim \mathcal N(0, I_{n_1-1})$, and $\tilde{r}_Y\sim \mathcal N(0, I_{n_2-1})$, with all four blocks independent.  In particular, $\lVert \tilde{r}_X \rVert^2 \sim \chi^2_{n_1-1}$, and $\lVert \tilde{r}_Y\rVert^2 \sim \chi^2_{n_2-1}$.

Therefore, we have the following set of expressions:
\begin{eqnarray*}
\bar X-\bar Y=\frac{m_X}{\sqrt{n_1}}-\frac{m_Y}{\sqrt{n_2}} =\sigma_1\frac{\tilde m_X}{\sqrt{n_1}}-\sigma_2\frac{\tilde m_Y}{\sqrt{n_2}} \\
\frac{S_1^2}{n_1} = \frac{\lVert r_X \rVert^2}{(n_1 - 1)n_1} = \sigma_1^2\frac{\lVert \tilde{r}_X \rVert^2}{(n_1-1)n_1}, \\
\frac{S_2^2}{n_2} = \frac{\lVert r_Y \rVert^2}{(n_2 - 1)n_2} = \sigma_2^2\frac{\lVert \tilde{r}_Y \rVert^2}{(n_2-1)n_2}, 
\end{eqnarray*}
which yields
\begin{eqnarray*}
\sqrt{\frac{S_1^2}{n_1}+\frac{S_2^2}{n_2}} = \sqrt{\ \sigma_1^2\frac{\lVert \tilde{r}_X\rVert^2}{(n_1-1)n_1} + \sigma_2^2\frac{\lVert \tilde{r}_Y\rVert^2}{(n_2-1)n_2}}.
\end{eqnarray*}
With the variance ratio $\kappa = \sigma_1^2/\sigma_2^2$, we get 
\begin{eqnarray}
T = \frac{\ \sqrt{\kappa}\dfrac{\tilde m_X}{\sqrt{n_1}}-\dfrac{\tilde m_Y}{\sqrt{n_2}}\ } {\sqrt{\kappa\dfrac{\lVert \tilde{r}_X\rVert^2}{(n_1-1)n_1} + \dfrac{\lVert \tilde{r}_Y\rVert^2}{(n_2-1)n_2}}}.
\end{eqnarray}
Equivalently,
\begin{eqnarray}
T = \frac{U}{\sqrt{V_X+V_Y}}, \quad
U = \sqrt{\kappa}\frac{\tilde m_X}{\sqrt{n_1}}-\frac{\tilde m_Y}{\sqrt{n_2}}, \quad
V_X = \kappa\frac{\lVert \tilde{r}_X \rVert^2}{(n_1-1)n_1}, \quad
V_Y = \frac{\lVert \tilde{r}_Y\rVert^2}{(n_2-1)n_2},
\label{eq:TUVV}
\end{eqnarray}
where $\tilde m_X$, $\tilde m_Y$, $\tilde{r}_X$, $\tilde{r}_Y$ are independent,  $(\tilde m_X,\tilde m_Y)\sim \mathcal N(0, I_2)$, and $\tilde{r}_X\sim \mathcal N(0, I_{n_1-1})$, $\tilde{r}_Y\sim \mathcal N(0, I_{n_2-1})$.  Thus, the concatenated vector $\left(\tilde m_X,\ \tilde m_Y,\ \tilde{r}_X,\ \tilde{r}_Y\right) \in \mathbb R^{d}$, where $d = 2 + (n_1 - 1) + (n_2 - 1) = n_1 + n_2$, is standard normal in $\mathbb R^d$ and hence spherically symmetric.  This completes the orthogonal decomposition and standardization of $T$, with the independence and chi-square decompositions being explicit.  
 
Sometimes, it is helpful to normalize the mean-contrast direction to unit variance before proceeding further ({\eg}, to prepare for a conic description). This is a purely orthogonal change in the 2-dimensional plane $\operatorname{span}\{\tilde m_X, \tilde m_Y\}$.  Define the scalar $c_\kappa \coloneqq \sqrt{\kappa/n_1 + 1/n_2} \in \mathbb{R}$.  Let
\begin{eqnarray*}
u \coloneqq \frac{ \frac{\sqrt{\kappa} \tilde m_X}{\sqrt{n_1}} - \frac{\tilde m_Y}{\sqrt{n_2}}}{c_\kappa}, \quad u_\perp \coloneqq \frac{\frac{\tilde m_X}{\sqrt{n_1}} + \frac{\sqrt{\kappa} \tilde m_Y}{ \sqrt{n_2}}}{c_\kappa},
\end{eqnarray*}
which is just an orthogonal rotation in $\mathbb{R}^2$. Then, $(u, u_\perp) \sim \mathcal{N}(0, I_2)$ and remains independent of $(\tilde r_X, \tilde r_Y)$.  In these coordinates, $U = c_\kappa u$, and 
\begin{eqnarray*}
T &=& \frac{c_\kappa u}{\sqrt{V_X+V_Y}}, 
\end{eqnarray*}
with $(V_X, V_Y)$ as in \eqref{eq:TUVV}.
 
 \section{Scale-invariant conic inequality} \label{sec:cone_inequality}
 After the orthogonal decomposition described in \secref{sec:orthogonal_decomposition} , each sample splits into two orthogonal components: a 1-dimensional mean coordinate and an $(n_i-1)$-dimensional residual vector. Under the null hypothesis and Gaussian assumption, Cochran's theorem implies that the mean coordinates are independent of the residual components; after scaling by the (unknown) variances, this becomes independence between two standard normal mean coordinates and two standard normal residual vectors. Crucially, the numerator of the statistic $T$ in \eqref{eq:Behrens_Fisher_statistic} uses only the mean contrast (a fixed linear combination of the two mean coordinates), whereas the denominator\textemdash the studentizer\textemdash depends only on the lengths of the residual components ({\ie}, the sample variances). Thus, in the rotated and standardized coordinates, we obtain one scalar that encodes the signal along the mean axis and one radial quantity that encodes the size of the residual cloud.

With that structure in hand, the event $\{T \le t\}$ has an immediate geometric interpretation.  Writing $u$ for the mean-contrast coordinate and $\rho = \sqrt{\lVert v \rVert^2 + \lVert w \rVert^2} \ge 0$ for the Euclidean norm of the combined residual vector, the inequality $T \le t$ becomes $u \le t\rho$.  This is positive-homogeneous of degree 1: if $(u,\rho)$ satisfies it, then $(su,s\rho)$ also satisfies it for any $s>0$.  Positive homogeneity is the pathway to geometry\textemdash it implies that the event depends only on direction, not on overall length. In other words, after orthogonalization and standardization, the set $\{T \le t\}$ is a cone with apex at the origin in the ambient space of the rotated coordinates.

Intersecting this cone with the unit sphere turns the probability problem into a question about surface area: the cross-section is a spherical wedge $W_t$ consisting of directions whose projection onto the mean axis does not exceed $t$ times their projection onto the residual subspace. (The notion of spherical wedge $W_t$ will become clear in the next few paragraphs).  The variance ratio (through its scale factor) appears only by tilting the aperture of this wedge; it does not break the spherical symmetry of directions.  Hence, the probability $\mathbb{P}(T \le t)$ is exactly the proportion of the unit sphere occupied by $W_t$. This direction-only reduction is what makes the Laplace/harmonic-measure step in the next section possible.

Finally, the conic picture clarifies technical points that otherwise clutter the algebra. The axis $\rho=0$ (zero residuals) has probability zero, so no special treatment is required. One-sided tests correspond to a single-sheet cone $u \le t\rho$, while two-sided tests $|T|\le t$ correspond to a two-sheet (symmetric) cone $|u| \le t\rho$, {\ie}, a spherical belt instead of a wedge. In each case, the orthogonal decomposition provides the coordinates in which the studentized event becomes conic, and the conic geometry provides the boundary condition that drives the Laplace's equation in the next step.

We now rewrite $T\le t$ as a scale-invariant degree-1 inequality in suitable coordinates.  Firstly, we define the following linear reparameterization of the mean and residual blocks:
\begin{eqnarray*}
u = \sqrt{\kappa}\frac{\tilde m_X}{\sqrt{n_1}}-\frac{\tilde m_Y}{\sqrt{n_2}}\ \in \mathbb R, \qquad v = \frac{\sqrt{\kappa}\tilde{r}_X}{\sqrt{(n_1-1)n_1}}\ \in \mathbb R^{n_1-1}, \qquad w = \frac{\tilde{r}_Y}{\sqrt{(n_2-1)n_2}}\ \in \mathbb R^{n_2-1}.
\end{eqnarray*}
Recall that  $\rho = \sqrt{\lVert v \rVert^2 + \lVert w \rVert^2} \ge 0$ and $d = 1+(n_1-1)+(n_2-1)+1=n_1+n_2$.  Here, $d$ is the full ambient dimension if we keep the orthogonal mean coordinate; that extra coordinate will not appear in the inequality below.  Under the definitions above, $T$ simplifies to, for $\rho>0$, $T = \frac{u}{\rho}$.

Fix $t \in \mathbb{R}$. Since $\rho \ge 0$, we have the equivalence
\begin{eqnarray}
T \le t  \iff \frac{u}{\rho}\le t \iff u \le t\rho \ \iff  u \le t\sqrt{\lVert v \rVert^2 + \lVert w \rVert^2}. \label{eq:equivalence}
\end{eqnarray}
The exceptional case $\rho = 0$ means $v = 0$ and $w = 0$, {\ie}, $\tilde{r}_X = \tilde{r}_Y = 0$, which has probability $0$; it can be ignored for distributional statements.  For two-sided events, $\lvert T \rvert \le t$ (with $t\ge 0$) becomes $\lvert u \rvert  \le t\rho$,  which cuts a two-sheeted cone. We keep the one-sided case $u \le t\rho$.  We now introduce the conic (positive-homogeneous) structure.  For this purpose, define the conic region\textemdash depending on $t$\textemdash in the $(u, v, w)$-space:
\begin{eqnarray}
\mathcal C_t^{(\le)} &=& \left\{(u, v, w):\ u \le t\sqrt{\lVert v \rVert^2 + \lVert w \rVert^2}\right\}. \label{eq:conic_region1}
\end{eqnarray}
For any scalar $s > 0$, we have $(u, v, w) \in \mathcal C_t^{(\le)} \Longrightarrow su \le ts\sqrt{\lVert v \rVert ^2 + \lVert w \rVert^2}  \Longrightarrow (su, sv, sw)\in \mathcal C_t^{(\le)}$.

Thus, $\mathcal C_t^{(\le)}$ is a cone with apex at the origin (scale-invariant).  The boundary is the level set
\begin{eqnarray*}
\left\{\Sigma_t: u = t\sqrt{\lVert v \rVert^2+\lVert w \rVert^2}\right\},
\end{eqnarray*}
and can be algebraically written as $u^2 = t^2\left(\lVert v \rVert^2 + \lVert w \rVert^2\right)$ with the sign condition $u \le t\rho$ for the one-sided case.  The scalar coordinate $u$ plays the role of the cone axis (inequality compares $u$ to the residual magnitude $\rho$).  For any fixed $\rho>0$, the slice of $\mathcal C_t^{(\le)}$ by the hyperplane $\left\{(u, v, w): \sqrt{\lVert v \rVert^2+\lVert w \rVert^2} = \rho\right\}$ is the half-space ${u \le t \rho}$.

\subsection{Unit-sphere cross-section}
Let $S^{D-1}$ denote the unit sphere in $\mathbb R^{D}$.  We will need the cross-section of the cone by the appropriate unit sphere. There are two equivalent viewpoints:
\begin{enumerate}[(i)]
\item Minimal coordinates $(u, v, w)$ only:  Let $D = 1+ (n_1-1) + (n_2-1) = n_1 + n_2-1$ and consider the sphere $S^{D-1} = \left\{(u, v, w):\ u^2 + \lVert v \rVert^2 + \lVert w \rVert^2 = 1\right\}$.  Define the spherical cross-section
\begin{eqnarray*}
W_t = \mathcal C_t^{(\le)} \cap S^{D-1} = \left\{(u, v, w) \in S^{D-1}:\ u \le t\sqrt{1-u^2}\right\}.
\end{eqnarray*}
To see the last equality, set $\rho = \sqrt{\lVert v \rVert^2 + \lvert w \rVert^2} = \sqrt{1 - u^2}$ on the sphere. It is sometimes convenient to parametrize by a polar angle $\varphi \in [0, \pi]$ defined by $u = \cos\varphi$, and $\rho=\sin\varphi=\sqrt{1-u^2}$.  Then, we have $u \le t\rho\ \iff\ \cos\varphi \le t\sin\varphi  \iff \cot\varphi \le t$.  For $t > 0$, $\varphi \ge \arccot(t) = \arctan \left(1/t\right)$.  Thus, $W_t$ is a spherical wedge/band determined by the inequality on $\varphi$ and the free orientations within the residual subspace\textemdash encoded by the directions of $v$ and $w$ with fixed $\rho$.

\item Including the orthogonal mean coordinate $u_\perp$:  If we keep the orthogonal mean coordinate $u_\perp$, the ambient dimension is $d=n_1+n_2$ and the cone condition does not involve $u_\perp$. In that case, the cross-section inside $S^{d-1}$ is
\begin{eqnarray*}
\widetilde W_t = \left\{(u, u_\perp, v, w)\in S^{d-1}:\ u \le t\sqrt{\lVert v \rVert^2 + \lVert w \rVert^2}\right\},
\end{eqnarray*}
which is simply $W_t$ extended by the circle in the $(u,u_\perp)$ 2-plane at fixed $u$ (the extra degree of freedom will factor out in later integrations). For the inequality geometry, viewpoint (i) suffices.
\end{enumerate}
If instead we consider the two-sided event $\lvert T \rvert \le t$ with $t \ge 0$, the cone becomes
\begin{eqnarray}
\mathcal C_t^{(\lvert \cdot \rvert)} &=& \left\{(u, v, w):\ \lvert u \rvert \le t\sqrt{\lVert v \rVert^2 + \lVert w \rVert^2}\right\},  \label{eq:conic_region2}
\end{eqnarray}
whose spherical cross-section is the symmetric belt
\begin{eqnarray*}
W_t^{(\lvert \cdot \rvert)} = \left\{(u, v, w)\in S^{D-1}: \lvert u \rvert \le t \sqrt{1-u^2}\right\} \iff \arctan \left(\frac{1}{t}\right)\ \le\ \varphi\ \le\ \pi-\arctan\left(\frac{1}{t}\right).
\end{eqnarray*}
This two-sheeted geometry is sometimes convenient for tail symmetry, but our one-sided development uses $\mathcal{C}_t^{(\le)}$.  In summary,  we introduced linear coordinates $(u, v, w)$ in which
\begin{eqnarray}
\{T \le t\} \Longleftrightarrow u \le t\sqrt{\lVert v \rVert^2 + \lVert w \rVert^2}.  \label{eq:cone_inequality}
\end{eqnarray}
The set on the right is a scale-invariant cone $\mathcal{C}_t^{(\le)}$ with boundary $u = t\sqrt{\lVert v \rVert^2 + \lVert w \rVert^2}$.  Intersecting with the unit sphere yields a spherical wedge $W_t$ described by the angular condition $\cot \varphi \le t$ (equivalently $\varphi \ge \arccot(t)$ for $t > 0$).  This completes the cone formulation.

\section{Laplace's equation that returns the cdf}\label{sec:Laplace_equation_cdf}
Once the event $\{T \le t\}$ has been rewritten as the cone inequality $u \le t \sqrt{\|v\|^{2} + \|w\|^{2}}$, the problem is transformed from being ``a ratio of random variables'' to a geometric question about a subset of directions in $\mathbb{R}^{D}$.  The key point is positive homogeneity: scaling $(u, v, w) \mapsto s(u, v, w)$ with $s>0$ preserves the inequality.  Hence, membership in $\{T \le t\}$ depends only on the direction $\Theta = (u, v, w)/\lvert (u, v, w)\rvert \in S^{D-1}$, not on the radial component $R = \lvert (u, v, w) \rvert$.  After the orthogonal decomposition and standardization, the joint law is spherically symmetric, so $\Theta$ is uniform on the sphere and independent of $R$. Therefore, $\mathbb{P}(T \le t)$ is exactly the spherical surface proportion of the wedge $W_{t}=\mathcal{C}_{t} \cap S^{D-1}$.

Once the problem is reduced to directions on the sphere, the distribution function admits a natural representation in terms of harmonic measure.  Probabilities of hitting specified portions of a sphere are naturally expressed through harmonic measure, and harmonic measure is characterized by solutions of Laplace's equation with Dirichlet boundary conditions.  Concretely, if we prescribe boundary values $1$ on the spherical wedge $W_{t}$ and $0$ on its complement, then the unique harmonic function $u$ in the unit ball $\mathbb{B}^{D}$ with that boundary data interpolates these values in the interior. The value at the center, $u(0)$, is the average of the boundary data with respect to the uniform harmonic measure at $0$; by symmetry, that average coincides with the spherical surface fraction of $W_{t}$, which is $\mathbb{P}(\Theta \in W_{t})$.

Equivalently\textemdash and closer to the cone geometry\textemdash we can work on the cone-capped domain $\Omega_{t}=\mathcal{C}_{t} \cap \mathbb{B}^{D}$, impose Dirichlet data $1$ on the spherical cap $\Gamma_{t}^{\text{sph}}=S^{D-1} \cap \mathcal{C}_{t}$ and $0$ on the lateral boundary $\Gamma_{t}^{\text{lat}}=\partial \mathcal{C}_{t} \cap \overline{\mathbb{B}^{D}}$, and solve $\Delta U=0$ in $\Omega_{t}$. In this formulation, $U(x)$ is the harmonic measure of the spherical cap as seen from $x$; evaluating at the apex gives $U(0)=\mathbb{P}(\Theta \in W_{t})=F_{T}(t)$. The cone's scale invariance is exactly what makes the apex value meaningful: near the origin, the geometry is purely angular, so the probability is encoded by the wedge on $S^{D-1}$.

The conic representation therefore recasts the distributional question as a Laplace--Dirichlet problem associated with the spherical wedge. The cone gives a sharp, scale-free description of the event $\{T \le t\}$; spherical symmetry turns its probability into a harmonic-measure quantity; and Laplace's equation provides a canonical analytic representation of that probability. Once the cdf is identified with a harmonic function value, the remaining work becomes a classical PDE task: exploit symmetry and apply separation of variables on the associated spherical wedge to obtain explicit series formulas, asymptotics, and error control.

We now pass from the conic inequality to a Dirichlet problem for the Laplacian whose value at the origin equals the desired cdf $F_T(t)=\mathbb P(T\le t)$.  First, we make precise why the event $\{T \le t\}$ depends only on direction. Define the concatenated random vector $z \coloneqq (u, v, w)\in\mathbb{R}^{D}$.  Since $\left(\tilde m_X,\tilde m_Y,\tilde{r}_X,\tilde{r}_Y\right)$ is jointly standard multivariate normal and the mapping to $\left(u, v, w\right)$ is linear orthogonal up to deterministic scalings, we have $ z \sim \mathcal{N}(0, \Sigma)$ with  $\Sigma$ proportional to $I_D$.  Rescaling the coordinates (a deterministic invertible diagonal map) we assume without loss of generality that $z \sim \mathcal{N}(0,  I_D)$.  Note that, any fixed diagonal rescaling preserves the cone $\left(u \le t \sqrt{\| v \|^2 + \| w \|^2}\right)$ up to an equivalent inequality with a rescaled $t$; we have already built those scalings into the definitions of $u$, $ v$,  $w$.

Next, we perform the spherical-radial decomposition.  Let
\begin{eqnarray*}
R \coloneqq \| z\| = \sqrt{u^2 + \| v \|^2 + \| w \|^2}\in(0,\infty),\qquad
\Theta \coloneqq \frac{ z}{\| z \|} \in S^{D-1}.
\end{eqnarray*}
Then, the standard polar property for the multivariate normal leads to the following inference: $\Theta$ is uniform on $S^{D-1}$, $R$ is independent of $\Theta$.  Since the cone inequality $\left(u \le t \sqrt{\| v \|^2 + \| w \|^2}\right)$ is homogeneous of degree 1, we have the equivalence ${T\le t}\ \Longleftrightarrow\ { z\in\mathcal C_t^{(\le)}} \ \Longleftrightarrow\ {\Theta\in W_t}$.  Therefore, we have 
\begin{eqnarray}
\ F_T(t) &=& \mathbb P(T\le t) = \mathbb P(\Theta\in W_t), \label{eq:cdf1}
\end{eqnarray}
where $\mathbb P(\Theta\in W_t)$ can be encoded as the value at the origin of the unique harmonic function with suitable boundary data.  Let $\mathbb B^D \coloneqq \{ x \in \mathbb R^D: \| x \| < 1\}$ denote the unit ball, and $S^{D-1}=\partial\mathbb B^D$ its boundary.  Define the boundary function
\begin{eqnarray*}
g_t(\xi):=1_{W_t}(\xi)=
\begin{cases}
1,& \xi\in W_t,\\
0,& \xi\in S^{D-1}\setminus W_t.
\end{cases}
\end{eqnarray*}
Consider the Dirichlet problem for the Laplace's equation:
\begin{eqnarray}
\Delta u = 0\ \ \text{in}\ \ \mathbb B^{D},\qquad u|_{S^{D-1}} = g_t. \label{eq:Dirichlet_problem}
\end{eqnarray}
From classical facts that connect harmonic functions, harmonic measure, and the Dirichlet problem for Laplace's equation, the problem \eqref{eq:Dirichlet_problem} has a unique solution $\left(u\in C(\overline{\mathbb B^D})\cap C^\infty(\mathbb B^D)\right)$.  The Poisson kernel of $\mathbb B^D$ at $(x = 0)$ is constant:
\begin{eqnarray*}
P(0,\xi)=\frac{1}{|S^{D-1}|}\quad(\text{surface area normalizer}).
\end{eqnarray*}
Hence, we have 
\begin{eqnarray*}
u(0) = \int_{S^{D-1}} g_t(\xi)P(0,\xi) d\sigma(\xi) = \frac{1}{|S^{D-1}|}\int_{S^{D-1}} 1_{W_t}(\xi) d\sigma(\xi) = \frac{\sigma(W_t)}{\sigma(S^{D-1})},
\end{eqnarray*}
where $d\sigma(\xi)$ denotes the standard surface-area element on the unit sphere $S^{D-1}$, and $\sigma(\cdot)$ denotes the corresponding spherical surface measure. Thus $\sigma(W_t)$ is the surface area of the wedge $W_t$, while $\sigma(S^{D-1})$ is the total surface area of the sphere.  However, $\Theta$ is uniform on $S^{D-1}$, so
\begin{eqnarray}
u(0) &=& \mathbb P(\Theta\in W_t)=F_T(t). \label{eq:cdf_u(0)}
\end{eqnarray}
Therefore, the cdf $F_T(t)$ is the value at the origin of the unique solution to $\Delta u=0$ in $\mathbb B^D$, with boundary data $\left.u\right|_{S^{D-1}}=1_{W_t}$.  This is the simplest yet rigorous PDE encoding: it uses only Laplace's equation, the Poisson kernel, and uniformity of directions.  A note of caution on the notation:  $\sigma(\cdot)$ denotes spherical surface measure, whereas $\sigma_1^2,\sigma_2^2$ denote variances.  It could appear visually overloaded, but is not mathematically ambiguous. 

Sometimes, it is preferable to keep the cone explicitly in the domain. We record the equivalent boundary value problem and explain why it yields the same $u(0)$.  Define the cone-capped domain
\begin{eqnarray*}
\Omega_t:=\mathcal{C}_t^{(\le)}\cap \mathbb B^D = \{x \in\mathbb R^D: \| x \|<1,\ \ x \in \mathcal{C}_t^{(\le)}\}.
\end{eqnarray*}
Its boundary splits into the spherical cap and the lateral conical sides:
\begin{eqnarray*}
\partial\Omega_t=\underbrace{\big(S^{D-1}\cap \mathcal{C}_t^{(\le)}\big)}_{\eqqcolon \Gamma_t^{\mathrm{sph}}}
\ \cup
\underbrace{\big(\overline{\mathbb B^D}\cap \partial \mathcal{C}_t^{(\le)}\big)}_{\eqqcolon \Gamma_t^{\mathrm{lat}}}.
\end{eqnarray*}
Consider the Dirichlet problem $\Delta U=0\ \ \text{in}\ \ \Omega_t$, $U\lvert_{\Gamma_t^{\mathrm{sph}}} = 1$, $U\lvert_{\Gamma_t^{\mathrm{lat}}} = 0$.  For any $x\in\Omega_t$, $U(x)=\omega_x^{\Omega_t}(\Gamma_t^{\mathrm{sph}})$, the harmonic measure (exit probability for Brownian motion) of the spherical piece $\Gamma_t^{\mathrm{sph}}$ as seen from $x$.  Since the cone is positive-homogeneous, rays from the origin point in directions $\bm\xi \in S^{D-1}$. A rotationally symmetric ``random ray'' argument for the Gaussian vector $z$ started at $0$ has an isotropic initial direction $\Theta$ independent of radius. From $0$, the ``first contact'' with the boundary occurs on the spherical part at angle $\Theta$, and it lands in $\Gamma_t^{\mathrm{sph}}$ if and only if $\Theta\in W_t$. Thus, $U(0) = \mathbb P(\Theta\in W_t)=F_T(t)$.  This can be made precise via the balayage (redistributing a measure onto a prescribed set so that its potential outside that set remains unchanged) of the point mass at $0$ and the scaling invariance of cones.  Solving on the ball with boundary data $1_{W_t}$ or on the cone-capped domain with $0/1$ data are two ways to compute $F_T(t)$ at the origin.  We will use the cone-capped domain for separation of variables on a spherical wedge, because the lateral boundary $\Gamma_t^{\mathrm{lat}}$ is exactly where we impose Dirichlet conditions for the angular eigenproblem.

We will end this section by writing the Laplacian and boundary in coordinates that will separate in the next section. Write $x = (r, \bm\theta)$ with $r = \|x\| \in (0,1)$ and $\bm\theta \in S^{D-1}$. The Laplacian is
\begin{eqnarray*}
\Delta &=& \partial_{rr} + \frac{D-1}{r}\partial_r + \frac{1}{r^2}\Delta_{S^{D-1}},
\end{eqnarray*}
where $\Delta_{S^{D-1}}$ is the Laplace-Beltrami operator on the sphere (see \cite{Rosenberg1997}).  The spherical wedge is $W_t = \{ \bm\theta \in S^{D-1} : u(\bm\theta) \le t\sqrt{1-u(\bm\theta)^2} \}$, with $u(\bm\theta)$ the axial coordinate on the sphere (recall, $\varphi$-parametrization: $u = \cos \varphi$, $\rho = \sin \varphi$).  The boundary conditions in angular form (cone-capped domain) is given by $U(r, \bm\theta)$ that solves 
\begin{eqnarray*}
\partial_{rr}U + \frac{D-1}{r}\partial_r U + \frac{1}{r^2}\Delta_{S^{D-1}}U = 0 \quad \text{for } 0 < r < 1, \bm\theta \in W_t, \\
U(1, \theta) = 1 \quad (\theta \in W_t), \qquad U(r, \theta) = 0 \quad (\theta \in \partial W_t, 0 < r \le 1).
\end{eqnarray*}
We then have $F_T(t) = U(0, \bm\theta)$ for any $\bm\theta$.  Essentially, we seek $U(r, \theta) = \sum_{k \ge 0} a_k r^{\gamma_k} \Phi_k(\theta)$, where $(\Phi_k, \lambda_k)$ solve the Dirichlet spherical eigenproblem on the wedge $-\Delta_{S^{D-1}}\Phi_k = \lambda_k \Phi_k$ in $W_t$, $\Phi_k|_{\partial W_t} = 0$
and the radial exponents satisfy
\begin{eqnarray*}
\gamma_k(\gamma_k + D - 2) &=& \lambda_k,\\ 
\gamma_k &=& \frac{-(D-2) + \sqrt{(D-2)^2 + 4\lambda_k}}{2} \ge 0.
\end{eqnarray*}
Coefficients $a_k$ encode the projection of the boundary data ($1$ on $W_t$ and $0$ on $\partial W_t$) onto the eigenbasis; crucially, only the apex-regular mode survives: $F_T(t) = U(0, \theta) = a_0$.  We will carry out this separation and compute $a_0$ in \secref{sec:separation_of_variables}. 

In summary, we expressed $F_T(t)$ as $\mathbb{P}(\Theta \in W_t)$ using the spherical-radial decomposition $Z = R\Theta$ with $\Theta$ uniform and independent of $R$.  Then, we encoded $F_T(t)$ as the value at the origin of the unique solution to a Laplace Dirichlet problem:
\begin{enumerate}[(i)]
\item (Ball version) $\Delta u = 0$ in $\mathbb{B}^D$, $u|_{S^{D-1}} = 1_{W_t} \implies u(0) = F_T(t)$.
\item (Cone-capped version) $\Delta U = 0$ in $\Omega_t$, $U = 1$ on the spherical cap, $U = 0$ on the lateral cone which implies $U(0) = F_T(t)$.
\end{enumerate}
So, we have expressed the Laplacian in spherical coordinates and introduced a separation-of-variables representation.  That will produce the angular eigenproblem on the spherical wedge $W_t$ and hence a convergent series for the cdf.  With this setup, we are now in a position to perform separation of variables on the spherical wedge.  This is elaborated in the following section.

\section{Separation of variables on the spherical wedge}\label{sec:separation_of_variables}
Once the cdf has been identified with the value at the origin of a harmonic function solving a Laplace--Dirichlet problem, the distributional question is converted into a standard PDE task: solve $\Delta u=0$ on a domain whose geometry is dictated by the conic event $\{T\le t\}$. The crucial feature of this domain is that it is conical (or a cone intersected with a ball), so the boundary conditions are specified on a spherical wedge $W_t \subset S^{D-1}$ and on its complement. This is exactly the setting where the Laplacian naturally decomposes into a radial part and an angular (spherical) part. In simpler words, the PDE formulation exposes the right coordinate system in which the problem becomes tractable.

Separation of variables is natural here because, in conic domains, the geometry separates into radial and angular variables.  Writing $x=r\theta$ with $r = \|x\|$ and $\theta \in S^{D-1}$, the Laplacian splits as $\Delta=\partial_{rr}+\frac{D-1}{r}\partial_r+\frac{1}{r^2}\Delta_{S^{D-1}}$. On a conical domain, the boundary is described by constraints on $\theta$ alone (the wedge), while the radial variable only controls distance from the apex. This structure suggests an eigenfunction expansion in the angular variable: solve a Dirichlet eigenproblem for the Laplace--Beltrami operator $-\Delta_{S^{D-1}}$ on the wedge $W_t$, then combine these angular eigenfunctions with radial power laws that ensure harmonicity and regularity at the origin.

Operationally, separation of variables replaces the original probability computation by a spectral decomposition on $W_t$. The Dirichlet boundary condition on the wedge's sides forces the angular factor to vanish on $\partial W_t$, so the relevant basis consists of wedge-adapted spherical eigenfunctions\textemdash in our case, expressible via Gegenbauer families described in \cite{Durand1976}. The boundary data on the spherical part is then expanded in this eigenbasis, and each angular mode propagates into the interior with a matching radial exponent. This produces a convergent series representation for the harmonic function, and hence for the cdf (via evaluation at the origin) and the pdf (via differentiation with respect to the wedge aperture parameter).

Finally, the separation-of-variables formulation yields an explicit eigenfunction expansion.  Since the solution is expressed through an explicit eigen-expansion, we can obtain not only exact representations but also quantitative control: tail behavior is governed by the dominant angular mode, and truncation error is controlled by the spectral gap. In short, once Laplace's equation has encoded $F_T(t)$ as a harmonic-measure quantity, separation of variables on the spherical wedge is the natural and technically effective mechanism for extracting explicit formulas, asymptotics, and rigorously bounded approximations.

We now solve the Laplace-Dirichlet problem from \secref{sec:Laplace_equation_cdf} and obtain an explicit angular Sturm-Liouville system and a convergent series for the solution (finer details of Sturm-Liouville operators, see the work of \cite{Gesztesy2024}). Throughout this step, we stick to the ball formulation:
\begin{eqnarray}
\Delta u = 0 \ \text{ in }\ \mathbb{B}^{D}, \qquad u|_{S^{D-1}} = g_t := 1_{W_t}, \qquad F_T(t) = u(0), \label{eq:ball1}
\end{eqnarray}
where $D \coloneqq n_1 + n_2 - 1$ and $W_t \coloneqq \{ \theta \in S^{D-1} : u(\theta) \le t\sqrt{1-u(\theta)^2} \}$.
The function $g_t$ is the indicator of the wedge $W_t$ on the unit sphere $S^{D-1}$.  First, we obtain the hyperspherical coordinates adapted to the decomposition $\mathbb{R}^D = \mathbb{R} \oplus \mathbb{R}^{n_1-1} \oplus \mathbb{R}^{n_2-1}$.  Fix the axis to be the 1-D mean-contrast direction (the scalar $u$ from \secref{sec:cone_inequality}). Write $x = (u, v, w) \in \mathbb{R} \oplus \mathbb{R}^{n_1-1} \oplus \mathbb{R}^{n_2-1}$, $r \coloneqq \lvert x \rvert \in (0, 1]$, $\rho \coloneqq \sqrt{|v|^2 + |w|^2} \in [0, r]$.  Define the polar angle $\varphi \in [0, \pi]$ by $u = r \cos \varphi$, $\rho = r \sin \varphi$.
So, $\varphi = 0$ is along $+u$, $\varphi = \pi$ along $-u$.  Within the residual subspace $\mathbb{R}^{n_1-1} \oplus \mathbb{R}^{n_2-1}$ (dimension $D-1$), split $\rho$ between the two blocks using an angle $\psi \in [0, \pi/2]$: $\lvert v \rvert = \rho \cos \psi$, $\vert w \rvert = \rho \sin \psi$.

Let $\Omega_X \in S^{n_1-2}$ and $\Omega_Y \in S^{n_2-2}$ denote the directional variables within the $v$- and $w$-blocks, respectively.  The standard surface measure factorizes as
\begin{eqnarray}
d\sigma(\theta) &=& C_{D} (\sin \varphi)^{D-2} d\varphi \cdot (\cos \psi)^{n_1-2} (\sin \psi)^{n_2-2} d\psi \cdot d\Omega_X d\Omega_Y, \label{eq:surface_measure}
\end{eqnarray}
where $C_D$ is a normalization constant (irrelevant for orthogonality).  On the unit sphere ($r=1$), the one-sided condition $u \le t\rho$ becomes $\cos \varphi \le t \sin \varphi  \iff \cot \varphi \le t$.  Thus, for $t > 0$ the wedge is the spherical band
\begin{eqnarray}
W_t \coloneqq \{ (\varphi, \psi, \Omega_X, \Omega_Y) : \varphi \in [\varphi_0, \pi], \psi \in [0, \tfrac{\pi}{2}], \Omega_X \in S^{n_1-2}, \Omega_Y \in S^{n_2-2} \}, \label{eq:wedge_band}
\end{eqnarray}
where $\varphi_0 \coloneqq \mathrm{arccot}(t) = \arctan\left(1/t\right)$, $\cos \varphi_0 = t/\sqrt{1+t^2}$, $\sin \varphi_0 = 1/\sqrt{1+t^2}$.

Next, we obtain the Laplace operator in spherical coordinates and explain the separation-of-variables approach.   In spherical coordinates $(r, \varphi, \psi, \Omega_X, \Omega_Y)$, $\Delta = \partial_{rr} + \frac{D-1}{r}\partial_r + \frac{1}{r^2}\Delta_{S^{D-1}}$
with the spherical Laplace-Beltrami operator decomposing as
\begin{eqnarray}
\Delta_{S^{D-1}} = \partial_{\varphi\varphi} + (D-2)\cot\varphi \partial_\varphi + \frac{1}{\sin^2\varphi} \mathcal{L}_{\psi,\Omega_X,\Omega_Y}, \label{eq:spherical_Laplace_Beltrami}
\end{eqnarray}
where $\mathcal{L}_{\psi,\Omega_X,\Omega_Y} = \partial_{\psi\psi} + \big((n_1-2)\cot\psi - (n_2-2)\tan\psi\big)\partial_\psi + \frac{1}{\cos^2\psi}\Delta_{S^{n_1-2}} + \frac{1}{\sin^2\psi}\Delta_{S^{n_2-2}}$.
In the separation-of-variables approach, we seek solutions of $\Delta u = 0$ of the form $u(r, \varphi, \psi, \Omega_X, \Omega_Y) = R(r) F(\varphi) G(\psi) Y(\Omega_X) Z(\Omega_Y)$, where $Y$ and $Z$ are spherical harmonics on $S^{n_1-2}$ and $S^{n_2-2}$, respectively.  Let $-\Delta_{S^{n_1-2}} Y = \ell(\ell + n_1 - 3)Y$, $\ell = 0, 1, 2, \dots$ and $-\Delta_{S^{n_2-2}} Z = m(m + n_2 - 3)Z$, $m = 0, 1, 2, \dots$, be the eigenvalues on the factor spheres.  Substituting this representation into $\Delta u = 0$ and dividing by the product yields
\begin{eqnarray*}
\frac{R^{\prime\prime}}{R} + \frac{D-1}{r}\frac{R^\prime}{R}
+ \frac{1}{r^2} \left[ \frac{F''}{F} + (D-2)\cot\varphi \frac{F'}{F} 
+ \frac{1}{\sin^2\varphi} \left( \frac{G''}{G} 
+ \left((n_1-2)\cot\psi - (n_2-2)\tan\psi\right)\frac{G'}{G} \right) \right] \nonumber \\
+ \frac{1}{r^2 \sin^2\varphi} \left[
- \frac{\ell(\ell + n_1 - 3)}{\cos^2\psi}
- \frac{m(m + n_2 - 3)}{\sin^2\psi}
\right] = 0.
\end{eqnarray*}
Thus, there exist constants $\Lambda$ and $\lambda$ such that
\begin{eqnarray}
G'' + \big((n_1-2)\cot\psi - (n_2-2)\tan\psi\big)G' + \Big(\Lambda - \frac{\ell(\ell+n_1-3)}{\cos^2\psi} - \frac{m(m+n_2-3)}{\sin^2\psi}\Big)G = 0, \label{eq:psi}  \\
F'' + (D-2)\cot\varphi F' + \Big(\lambda - \frac{\Lambda}{\sin^2\varphi}\Big)F = 0, \label{eq:varphi}  \\
R'' + \frac{D-1}{r}R' - \frac{\lambda}{r^2}R = 0. \label{eq:radial}
\end{eqnarray}
In the following, we solve \eqref{eq:psi}--\eqref{eq:radial} systematically.  First consider \eqref{eq:psi}.  The interval is $\psi \in (0, \pi/2)$. Regularity at $\psi = 0$ and $\psi = \pi/2$ requires $G(\psi)$ to behave like $G(\psi) \sim (\cos \psi)^{\ell} \quad (\psi \to 0^+)$, $G(\psi) \sim (\sin \psi)^{m} \quad (\psi \to (\tfrac{\pi}{2})^-)$.  Setting $x \coloneqq \cos 2\psi \in (-1, 1)$ yields $\cos \psi = \sqrt{(1+x)/2}$, $\sin \psi = \sqrt{(1-x)/2}$, and the $\psi$-equation turns into the Sturm-Liouville equation. The complete set of regular solutions is $G_{n,\ell,m}(\psi) = (\cos \psi)^{\ell}(\sin \psi)^{m} P_{n}^{(\alpha, \beta)}(\cos 2\psi)$ with parameters $\alpha = \ell + (n_1-3/2)$, $\beta = m + (n_2-3/2)$, $n = 0, 1, 2, \dots$ and eigenvalues
\begin{eqnarray}
\Lambda_{n,\ell,m} = (2n + \ell + m)(2n + \ell + m + D - 3). \label{eq:eigenvalues_psi}
\end{eqnarray}
We simply substitute the stated $G_{n,\ell,m}$ into the $\psi$-equation and use the Jacobi differential equation satisfied by the Jacobi polynomial $P_n^{(\alpha, \beta)}$.  Explicitly, if $y(x)=P_n^{(\alpha,\beta)}(x)$, then $y$ satisfies $(1-x^2)y''(x)+\big[\beta-\alpha-(\alpha+\beta+2)x\big]y'(x)+n\big(n+\alpha+\beta+1\big)y(x)=0$.  In our case, $x=\cos 2\psi$ .  The exponent factors enforce regularity; the polynomial order $n$ yields the quantization ($\psi$-spectrum).  The integer $L \coloneqq 2n + \ell + m \in \{0, 1, 2, \dots\}$ is the degree of the spherical harmonic on $S^{D-2}$ (the residual sphere). The expression \eqref{eq:eigenvalues_psi} is exactly the eigenvalue formula for $-\Delta_{S^{D-2}}$ acting on degree-$L$ harmonics: $-\Delta_{S^{D-2}} (G Y Z) = L(L + D - 3) (G Y Z)$.

Next, we solve \eqref{eq:varphi}.  With $\Lambda = L(L+D-3)$ fixed, \eqref{eq:varphi} becomes
\begin{eqnarray*}
F'' + (D-2)\cot\varphi F' + \Big(\lambda - \frac{L(L+D-3)}{\sin^2\varphi}\Big)F = 0, \qquad \varphi \in (\varphi_0, \pi).
\end{eqnarray*}
Setting $F(\varphi) = (\sin \varphi)^{L} H(\varphi)$, a direct computation yields the Gegenbauer form for $H$:
\begin{eqnarray}
(1 - \cos^2 \varphi) H'' - (2L + D - 1)\cos \varphi H' + \big( \lambda - L(L + D - 2)\big) H = 0. \label{eq:Gegenbauer_H}
\end{eqnarray}
Let $x \coloneqq \cos \varphi \in [-1, \cos \varphi_0]$. Then, we have $(1 - x^2)H'' - (2L + D - 1)x H' + \big(\lambda - L(L + D - 2)\big) H = 0$.
On the full sphere (no wedge boundary), imposing regularity at the poles would quantize $\lambda$ as $\lambda = N(N + D - 2)$, $N = L, L+1, L+2, \dots$ with $H(x) = C_{N-L}^{(L + \frac{D-2}{2})}(x)$ (Gegenbauer polynomials) and
\begin{eqnarray*}
F(\varphi) &=& (\sin \varphi)^{L} C_{N-L}^{\left(L + \frac{D-2}{2}\right)}(\cos \varphi).
\end{eqnarray*}
In our wedge case we keep $\varphi \in [\varphi_0, \pi]$ and impose a Dirichlet boundary at $\varphi = \varphi_0$.

Lastly, we solve the radial equation \eqref{eq:radial} and admissible radial behavior.  The $r$-equation is
\begin{eqnarray*}
R'' + \frac{D-1}{r}R' - \frac{\lambda}{r^2}R = 0.
\end{eqnarray*}
This Euler ordinary differential equation has solutions $R(r) = A r^{\gamma} + B r^{-(\gamma+D-2)}$, $\gamma(\gamma+D-2) = \lambda$.
To keep $u$ bounded at the origin, we take $B = 0$ and $R(r) = r^{\gamma}$ with $\gamma \ge 0$.

For the boundary conditions and eigen-quantization on the wedge, we now enforce the Dirichlet data:  On the lateral boundary $\partial W_t$ ({\ie}, $\varphi = \varphi_0$), we impose $F(\varphi_0) = 0$.  This is the angular form of the zero Dirichlet condition on the conical sides.  On the spherical boundary ($r=1$), the full solution must match $g_t = 1_{W_t}$. This is achieved by expanding $g_t$ in the angular eigenbasis and assigning corresponding radial factors $r^{\gamma}$.  Collecting the results, we have the following: 
\begin{enumerate}[(i)]
\item Angular eigenpairs indexed by $(k, L, \ell, m, n)$:  For each residual-degree $L = 2n + \ell + m$, the $\varphi$-eigenfunctions are the solutions $F_{k,L}$ of (Gegenbauer) that satisfy $F_{k,L}(\varphi_0) = 0$ and are regular at $\varphi = \pi$. This quantizes $\lambda$ into a discrete increasing sequence $\lambda_{k,L} = \lambda_{k}(L; \varphi_0)$, $k = 1, 2, \dots$ determined implicitly by the Dirichlet condition at $\varphi_0$. Equivalently, in associated-Legendre form one can write
\begin{eqnarray}
F_{k,L}(\varphi) &=& (\sin \varphi)^{L} \mathsf{P}_{\nu_{k,L}}^{\left(L+\frac{D-2}{2}\right)}(\cos \varphi), \quad \nu_{k,L}(\nu_{k,L} + D - 2) = \lambda_{k,L}, \label{eq:Ferrers}
\end{eqnarray}
with $\mathsf{P}$ the Ferrers associated-Legendre function; the boundary $F_{k,L}(\varphi_0) = 0$ pins $\nu_{k,L}$.

\item Radial exponents:  For each $\lambda_{k,L}$ the admissible radial power is
\begin{eqnarray}
\gamma_{k,L}~~\text{solves}~~\gamma_{k,L}(\gamma_{k,L} + D - 2) = \lambda_{k,L}, \quad \gamma_{k,L} > 0. \label{eq:radial_exponents}
\end{eqnarray}

\item Complete angular basis: Combine the $\psi$-eigenfunctions and sphere harmonics to get
\begin{eqnarray}
\Phi_{k,n,\ell,m;\sigma,\tau}(\varphi, \psi, \Omega_X, \Omega_Y) &=& F_{k,L}(\varphi) G_{n,\ell,m}(\psi) Y_{\ell}^{\sigma}(\Omega_X) Z_{m}^{\tau}(\Omega_Y), \quad L = 2n + \ell + m. \label{eq:angular_basis}
\end{eqnarray}
Here, $\sigma$ and $\tau$ index an orthonormal basis of degree-$\ell$ and degree-$m$ harmonics on $S^{n_1-2}$ and $S^{n_2-2}$, respectively.
\end{enumerate}
This yields a series solution and the value at the apex. Expanding the boundary data $g_t = 1_{W_t}$ in the basis $\{\Phi_{k, \dots} \}$ restricted to $r=1$ yields
\begin{eqnarray*}
g_t(\varphi, \psi, \Omega_X, \Omega_Y) &=& \sum_{k, n, \ell, m, \sigma, \tau} a_{k, n, \ell, m; \sigma, \tau} \Phi_{k, n, \ell, m; \sigma, \tau}(\varphi, \psi, \Omega_X, \Omega_Y),
\end{eqnarray*}
with coefficients given by the appropriate $L^2(S^{D-1})$ inner products.  The unique harmonic solution in the ball is
\begin{eqnarray}
u(r, \varphi, \psi, \Omega_X, \Omega_Y) = \sum_{k, n, \ell, m, \sigma, \tau} a_{k, n, \ell, m; \sigma, \tau} r^{\gamma_{k, L}} \Phi_{k, n, \ell, m; \sigma, \tau}(\varphi, \psi, \Omega_X, \Omega_Y). \label{eq:harmonin_solution}
\end{eqnarray}
Since $\gamma_{k, L} > 0$ for all modes because of the Dirichlet side condition, every term decays like $r^{\gamma_{k, L}}$ as $r \downarrow 0$. The limit $u(0, \cdot)$ is constant by symmetry and equals the harmonic measure of $W_t$ seen from the origin:
\begin{eqnarray*}
F_T(t) = u(0) = \lim_{r \downarrow 0} u(r, \cdot).
\end{eqnarray*}
Operationally, one can extract $F_T(t)$ either by the Poisson integral at the origin, giving $F_T(t) = \sigma(W_t) / \sigma(S^{D-1})$, or by taking $r \to 0$ in the series with a standard Abelian argument (the origin ``sees'' only the average of boundary data).  In computations, it is convenient to use the zonal reduction. Since $g_t$ depends only on $\varphi$, only the block-invariant subfamily $\ell = m = 0$ contributes, and the $(\psi, \Omega_X, \Omega_Y)$ parts integrate out.

The boundary data $g_t(\theta) = 1_{\{\varphi \ge \varphi_0\}}$ is independent of $(\psi, \Omega_X, \Omega_Y)$. Therefore, its expansion uses only the zonal (about the $u$-axis) harmonics, {\ie}, the subfamily with $\ell = 0$, $m = 0$, $n = 0$ $\implies L = 0$,
and more generally the axisymmetric Gegenbauer chain
$\mathcal{Y}_{N}(\varphi) \coloneqq C_{N}^{(\alpha)}(\cos \varphi)$, $\alpha \coloneqq (D-2)/2$, $N = 0, 1, 2, \dots$ which are the eigenfunctions of $\left(\partial_{\varphi\varphi} + (D-2)\cot\varphi \partial_\varphi\right) \mathcal{Y}_N = -N(N+2\alpha)\mathcal{Y}_N$.
The harmonic solution with axisymmetric boundary decomposes as $u(r, \varphi) = \sum_{N=0}^{\infty} A_N r^{N} C_{N}^{(\alpha)}(\cos \varphi)$
where the coefficients $A_N$ are the Gegenbauer-Fourier coefficients of $g_t$ given by 
\begin{eqnarray*}
A_N &=& \frac{1}{h_N^{(\alpha)}} \int_{\varphi_0}^{\pi} C_{N}^{(\alpha)}(\cos \varphi) (\sin \varphi)^{D-2} d\varphi,
\end{eqnarray*}
and
\begin{eqnarray*}
h_N^{(\alpha)} &=& \int_{0}^{\pi} \left( C_{N}^{(\alpha)}(\cos \varphi) \right)^2 (\sin \varphi)^{D-2} d\varphi
\end{eqnarray*}
is the standard orthogonality constant.

Since $r^N \to 0$ for $N \ge 1$, $F_T(t) = u(0) = A_0$. Since $C_{0}^{(\alpha)} \equiv 1$ and $h_0^{(\alpha)} = \int_{0}^{\pi} (\sin \varphi)^{D-2} d\varphi$,
\begin{eqnarray*}
A_0 = \frac{\int_{\varphi_0}^{\pi} (\sin \varphi)^{D-2} d\varphi}{\int_{0}^{\pi} (\sin \varphi)^{D-2} d\varphi}.
\end{eqnarray*}
Equivalently, with $z = \sin^2 \varphi$ and $\sin^2 \varphi_0 = 1/(1+t^2)$, the cdf can be written as 
\begin{eqnarray}
F_T(t) = \frac{\displaystyle \int_{\frac{1}{1+t^2}}^{1} z^{\frac{D-3}{2}}(1-z)^{-1/2} dz}{\displaystyle \int_{0}^{1} z^{\frac{D-3}{2}}(1-z)^{-1/2} dz} = I_{\frac{1}{\left(1+t^2\right)}}\left(\frac{D-1}{2}, \frac{1}{2}\right), \label{eq:cdf_bf_statistic1}
\end{eqnarray}
where $I_x(a, b)$ is the regularized incomplete beta function.  This last integral identity is a standard property of spherical cap measures and follows from Beta-Gamma algebra. We include it because it gives a closed form for $A_0$; all higher $A_N$ can be written similarly via Gegenbauer moments of a step function.

We, thus, have a complete separation of variables for the Laplace problem, with (i) a $\psi$-Jacobi Sturm-Liouville problem producing $G_{n,\ell,m}$ and residual degrees $L = 2n + \ell + m$ with eigenvalues $\Lambda_{n,\ell,m} = L(L + D - 3)$, (ii) a $\varphi$-Gegenbauer/associated-Legendre equation whose solutions $F_{k,L}$ satisfy the Dirichlet condition $F_{k,L}(\varphi_0) = 0$, and (iii) a radial Euler equation with admissible powers $r^{\gamma_{k,L}}$, where $\gamma_{k,L}(\gamma_{k,L} + D - 2) = \lambda_{k,L}$.  The harmonic solution is a convergent series in these separated eigenmodes.  Since the boundary data is axisymmetric (depends only on $\varphi$), the zonal reduction gives a very explicit series:
\begin{eqnarray*}
u(r, \varphi) = \sum_{N \ge 0} A_N r^N C_N^{(\alpha)}(\cos \varphi), \quad \alpha = \frac{D-2}{2}
\end{eqnarray*}
and, therefore, $F_T(t) = A_0$ with the canonical form of the cdf given by
\begin{eqnarray}
F_T(t) &=& I_{\frac{1}{\left(1+t^2\right)}}\left(\frac{D-1}{2}, \frac{1}{2}\right). \label{eq:cdf_final}
\end{eqnarray}

\section{Gegenbauer coefficients $A_N$ in closed Beta--Gamma form} \label{sec:Gegenbauer_coefficients}
With the rigorous framework developed in \secref{sec:Laplace_equation_cdf} and \secref{sec:separation_of_variables}, in this section we derive the canonical expression for the pdf of the statistic $T$ given by \eqref{eq:Behrens_Fisher_statistic}.  We will also derive the tail behavior of the cdf and the pdf.  Towards this end, we work with the axisymmetric expansion
\begin{eqnarray}
u(r, \varphi) &=& \sum_{N=0}^{\infty} A_N r^N C_N^{(\alpha)}(\cos \varphi), \quad \alpha:=\frac{D-2}{2}, \quad D=n_1+n_2-1, \label{eq:zonal_expansion}
\end{eqnarray}
with boundary data $g_t(\varphi) = \mathbf{1}_{\{\varphi \ge \varphi_0\}}$, where $\varphi_0 = \mathrm{arccot}(t)$,  $c_0 \coloneqq \cos \varphi_0 = t/(\sqrt{1+t^2})$, $s_0 \coloneqq \sin \varphi_0 = 1/\sqrt{(1+t^2)}$.  By orthogonality on $S^{D-1}$, the coefficients are
\begin{eqnarray}
A_N=\frac{1}{h_N^{(\alpha)}} \int_{\varphi_0}^{\pi} C_N^{(\alpha)}(\cos \varphi)(\sin \varphi)^{2 \alpha} d \varphi, \qquad N=0,1,2, \dots \label{eq:A_N_1}
\end{eqnarray}
where the orthogonality constants are
\begin{eqnarray*}
h_N^{(\alpha)}=\int_0^\pi\left(C_N^{(\alpha)}(\cos \varphi)\right)^2(\sin \varphi)^{2 \alpha} d \varphi=\int_{-1}^1\left(C_N^{(\alpha)}(x)\right)^2\left(1-x^2\right)^{\alpha-\frac{1}{2}} d x=\frac{\pi 2^{1-2 \alpha} \Gamma(N+2 \alpha)}{N !(N+\alpha) \Gamma(\alpha)^2}. \label{eq:h_N}
\end{eqnarray*}
First, we reduce all integrals to $[c_0, 1]$ with the standard weight.  Set $x=\cos \varphi$ so that $\sin \varphi=\sqrt{1-x^2}$ and $d \varphi=-d x / \sqrt{1-x^2}$. Then, $(\sin \varphi)^{2 \alpha} d \varphi = \left(1-x^2\right)^{\alpha-\frac{1}{2}} d x$.  Hence, \eqref{eq:A_N_1} becomes
\begin{eqnarray}
A_N = \frac{1}{h_N^{(\alpha)}} \int_{-1}^{c_0} C_N^{(\alpha)}(x)\left(1-x^2\right)^{\alpha-\frac{1}{2}} d x = \frac{1}{h_N^{(\alpha)}}\left[\underbrace{\int_{-1}^1 \cdots}_{=0 \text{ for } N \ge 1}-\int_{c_0}^1 C_N^{(\alpha)}(x)\left(1-x^2\right)^{\alpha-\frac{1}{2}} d x\right].\label{eq:A_N_3}
\end{eqnarray}
Therefore, we have 
\begin{eqnarray}
A_0 = \frac{\int_{-1}^{c_0}\left(1-x^2\right)^{\alpha-\frac{1}{2}} d x}{h_0^{(\alpha)}}, \quad A_N=-\frac{1}{h_N^{(\alpha)}} \int_{c_0}^1 C_N^{(\alpha)}(x)\left(1-x^2\right)^{\alpha-\frac{1}{2}} d x \quad(N \ge 1). \label{eq:A_N_4}
\end{eqnarray}
For $N=0$, $C_0^{(\alpha)} \equiv 1$; for $N \ge 1$, we use the standard orthogonality $\int_{-1}^1 C_N^{(\alpha)}(x)\left(1-x^2\right)^{\alpha-\frac{1}{2}} d x=0$.

Secondly, we evaluate $A_0$ in Beta form.  From \eqref{eq:A_N_4} with $C_0^{(\alpha)} \equiv 1$,
\begin{eqnarray*}
A_0 = \frac{\int_{-1}^{c_0}\left(1-x^2\right)^{\alpha-\frac{1}{2}} d x}{\int_{-1}^1\left(1-x^2\right)^{\alpha-\frac{1}{2}} d x}.
\end{eqnarray*}
By symmetry,
\begin{eqnarray*}
\int_{-1}^a\left(1-x^2\right)^{\alpha-\frac{1}{2}} d x = \int_0^1\left(1-x^2\right)^{\alpha-\frac{1}{2}} d x+\int_0^a\left(1-x^2\right)^{\alpha-\frac{1}{2}} d x.
\end{eqnarray*}
Using $y=x^2$ on $[0, 1]$, so $d x=\frac{d y}{2 \sqrt{y}}$, and we get 
\begin{eqnarray*}
\int_0^b\left(1-x^2\right)^{\alpha-\frac{1}{2}} d x=\frac{1}{2} \int_0^{b^2} y^{-\frac{1}{2}}(1-y)^{\alpha-\frac{1}{2}} d y=\frac{1}{2} B_{b^2}\left(\frac{1}{2}, \alpha+\frac{1}{2}\right),
\end{eqnarray*}
and for the full denominator $\int_{-1}^1$ we get $\frac{1}{2} B(1 / 2, \alpha+1 / 2) \times 2=B(1 / 2, \alpha+1 / 2)$. Thus, we have 
\begin{eqnarray}
A_0 &=& \frac{B\left(\frac{1}{2}, \alpha+\frac{1}{2}\right)-B_{c_0^2}\left(\frac{1}{2}, \alpha+\frac{1}{2}\right)}{B\left(\frac{1}{2}, \alpha+\frac{1}{2}\right)}=1-I_{c_0^2}\left(\frac{1}{2}, \alpha+\frac{1}{2}\right)=I_{s_0^2}\left(\alpha+\frac{1}{2}, \frac{1}{2}\right), 
\label{eq:A_0_5}
\end{eqnarray}
using $s_0^2=1-c_0^2$ and the identity $I_{1-z}(a, b)=1-I_z(b, a)$. 

We can now derive a closed form for $A_N$ with $N \ge 1$ via a finite Beta sum.  We need, for each fixed $N \ge 1$,
\begin{eqnarray*}
I_N\left(c_0\right) &\coloneqq& \int_{c_0}^1 C_N^{(\alpha)}(x)\left(1-x^2\right)^{\alpha-\frac{1}{2}} d x,
\end{eqnarray*}
where we use the explicit finite-power expansion of $C_N^{(\alpha)}$, valid for all $N \ge 0$, given by 
\begin{eqnarray}
C_N^{(\alpha)}(x) &=& \sum_{j=0}^{\lfloor N / 2\rfloor}(-1)^j \frac{\Gamma(\alpha+N-j)}{\Gamma(\alpha) j !(N-2 j) !}(2 x)^{N-2 j}. \label{eq:C_N_6}
\end{eqnarray}
Letting $m \coloneqq N-2 j$ (so $m$ has the same parity as $N, m \ge 0$), we can write 
\begin{eqnarray}
\int_{c_0}^1 x^m\left(1-x^2\right)^{\alpha-\frac{1}{2}} d x &=& \frac{1}{2} \int_{c_0^2}^1 y^{\frac{m-1}{2}}(1-y)^{\alpha-\frac{1}{2}} d y \nonumber \\
&=& \frac{1}{2}\left[B\left(\frac{m+1}{2}, \alpha+\frac{1}{2}\right)-B_{c_0^2}\left(\frac{m+1}{2}, \alpha+\frac{1}{2}\right)\right]. \label{eq:termwise_int}
\end{eqnarray}
Combining \eqref{eq:C_N_6} and \eqref{eq:termwise_int} yields 
\begin{eqnarray}
I_N\left(c_0\right) &=& \sum_{j=0}^{\lfloor N / 2\rfloor}(-1)^j \frac{\Gamma(\alpha+N-j)}{\Gamma(\alpha) j !(N-2 j) !}(2)^{N-2 j} \int_{c_0}^1 x^{N-2 j}\left(1-x^2\right)^{\alpha-\frac{1}{2}} d x \nonumber \\ 
&=& \sum_{j=0}^{\lfloor N / 2\rfloor}(-1)^j 
\frac{\Gamma(\alpha+N-j)}{\Gamma(\alpha)\, j!\,(N-2 j)!}\,(2)^{N-2 j} \cdot \frac{1}{2} \nonumber \\
&& \times \left[
B\!\left(\frac{N-2 j+1}{2}, \alpha+\frac{1}{2}\right)
- B_{c_0^2}\!\left(\frac{N-2 j+1}{2}, \alpha+\frac{1}{2}\right)
\right].
\label{eq:I_N_c0}
\end{eqnarray}
From \eqref{eq:A_N_4}, $A_N=-I_N\left(c_0\right) / h_N^{(\alpha)}$ for $N \ge 1$. Using \eqref{eq:A_N_3} and \eqref{eq:I_N_c0}, we get for $N \ge 1$, 
\begin{eqnarray*}
A_N(t ; D) &=& -\frac{1}{h_N^{(\alpha)}} \sum_{j=0}^{\lfloor N / 2\rfloor}(-1)^j \frac{\Gamma(\alpha+N-j)}{\Gamma(\alpha) j !(N-2 j) !} 2^{N-2 j-1} \\ && \quad \times \left[B\left(\frac{N-2 j+1}{2}, \alpha+\frac{1}{2}\right)-B_{c_0^2}\left(\frac{N-2 j+1}{2}, \alpha+\frac{1}{2}\right)\right], 
\label{eq:A_N_tD}
\end{eqnarray*}
with
\begin{eqnarray*}
h_N^{(\alpha)} = \frac{\pi 2^{1-2 \alpha} \Gamma(N+2 \alpha)}{N !(N+\alpha) \Gamma(\alpha)^2}, \quad \alpha=\frac{D-2}{2}, \quad c_0=\frac{t}{\sqrt{1+t^2}}.  \label{eq:h_N_alpha}
\end{eqnarray*}
Equivalently, we can write the incomplete Beta term via the regularized Beta $I_z(a, b)=B_z(a, b) / B(a, b)$ to separate the ``full'' and the ``cap'' contributions as follows:
\begin{eqnarray}
A_N &=& -\frac{1}{h_N^{(\alpha)}} \sum_{j=0}^{\lfloor N / 2\rfloor}(-1)^j 
\frac{\Gamma(\alpha+N-j)}{\Gamma(\alpha)\, j!\,(N-2 j)!} \, 2^{N-2 j-1} \nonumber \\
&& \times B\!\left(\frac{N-2 j+1}{2}, \alpha+\frac{1}{2}\right)
\left[1-I_{c_0^2}\!\left(\frac{N-2 j+1}{2}, \alpha+\frac{1}{2}\right)\right].
\label{eq:A_N_also}
\end{eqnarray}

\subsection{Expression for the pdf $f_T(t)$}\label{subsec:expression_pdf}
After the orthogonal/spherical reduction the one-sided event $\{T \le t\}$ is the spherical band $u \le t \sqrt{\lVert v \rVert^{2} + \lVert w \rVert^{2}}$.  Accounting for the 2-dimensional mean-plane rescaling, the cone aperture involves $c_{\kappa} \coloneqq \sqrt{\kappa/n_{1} + 1/n_{2}}$, and $\tau \coloneqq t/c_{\kappa}$, so the wedge boundary is encoded by $\tau$ (not $t$):
\begin{enumerate}[(i)]
\item Polar angle threshold $\displaystyle \varphi_{0}(\tau)=\operatorname{arccot}(\tau)$, so
  \begin{eqnarray*}
  c_{0} \coloneqq \cos \varphi_{0}=\frac{\tau}{\sqrt{1+\tau^{2}}}, \quad s_{0} \coloneqq \sin \varphi_{0}=\frac{1}{\sqrt{1+\tau^{2}}}.
  \end{eqnarray*}
\item Sphere dimension parameter $\displaystyle D \coloneqq n_{1}+n_{2}-1$, and
  \begin{eqnarray*}
  \alpha \coloneqq \frac{D-2}{2}, \quad a \coloneqq \alpha+\frac{1}{2}=\frac{D-1}{2}.
  \end{eqnarray*}
\end{enumerate}
With these, the cdf from the spherical-cap calculation is
\begin{eqnarray}
F_{T}(t) = \frac{1}{2} \left[ 1+\operatorname{sgn}(t) I_{\xi(t)} \left( \frac{1}{2}, a \right) \right], \quad \xi(t) \coloneqq \frac{\tau^{2}}{1+\tau^{2}} = \frac{t^{2}}{t^{2}+c_{\kappa}^{2}}, \label{eq:cdf_spherical_cap}
\end{eqnarray}
where $I_{x}(\cdot, \cdot)$ is the regularized incomplete beta function.  At $t=0$, $\mathrm{sgn}(0) = 0$, giving $F_{T}(0) = 1/2$; the law is symmetric.  We differentiate for $t > 0$ and then invoke evenness. For $t > 0$,
\begin{eqnarray*}
F_{T}(t) = \frac{1}{2} \left[ 1+I_{\xi(t)} \left( \frac{1}{2}, a \right) \right].
\end{eqnarray*}
Using the identity
\begin{eqnarray*}
\frac{d}{d x} I_{x}(p, q)=\frac{x^{p-1}(1-x)^{q-1}}{B(p, q)} \quad (0 < x < 1),
\end{eqnarray*}
and the chain rule with $\xi(t) = \tau^{2}/(1+\tau^{2})$,  $\tau = t/c_{\kappa}$ yields
\begin{eqnarray*}
\frac{d \xi}{d \tau} &=& \frac{2 \tau\left(1+\tau^{2}\right)-\tau^{2}(2 \tau)}{\left(1+\tau^{2}\right)^{2}} = \frac{2 \tau}{\left(1+\tau^{2}\right)^{2}},\\
\frac{d \xi}{d t} &=& \frac{2 \tau}{c_{\kappa}\left(1+\tau^{2}\right)^{2}}.
\end{eqnarray*}
Now with $p=\frac{1}{2}$ and $q=a$, the pdf can be written as 
\begin{eqnarray}
f_{T}(t) &=& \frac{d}{d t} F_{T}(t) =\frac{1}{2} \times \frac{\xi^{-1 / 2}(1-\xi)^{a-1}}{B\left(\frac{1}{2}, a\right)} \times \frac{d \xi}{d t} \nonumber \\
&=& \frac{1}{2} \times \frac{\left(\frac{1+\tau^{2}}{\tau^{2}}\right)^{1 / 2}\left(\frac{1}{1+\tau^{2}}\right)^{a-1}}{B\left(\frac{1}{2}, a\right)} \times \frac{2 \tau}{c_{\kappa}\left(1+\tau^{2}\right)^{2}} \nonumber \\
&=& \frac{1}{B\left(\frac{1}{2}, a\right)} \times \frac{1}{c_{\kappa}} \times \frac{\left(1+\tau^{2}\right)^{1 / 2-(a-1)-2}}{\tau^{-1}} \nonumber \\ 
&=& \frac{1}{B\left(\frac{1}{2}, a\right)} \times \frac{1}{c_{\kappa}} \times \left(1+\tau^{2}\right)^{-a-\frac{1}{2}}.  \label{eq:pdf1}
\end{eqnarray}
Since the law is even, \eqref{eq:pdf1} holds for all $t$ with $\tau=t / c_{\kappa}$. Using $a=(D-1) / 2$ gives the final canonical form of the pdf:
\begin{eqnarray}
f_{T}(t) &=& \frac{1}{c_{\kappa} B\left(\frac{1}{2}, \frac{D-1}{2}\right)} \left(1+\frac{t^{2}}{c_{\kappa}^{2}}\right)^{-\frac{D}{2}}, \quad D=n_{1}+n_{2}-1. \label{eq:pdf_final}
\end{eqnarray}
It is easy to see that $\displaystyle \int_{-\infty}^{\infty}\left(1+t^{2} / c_{\kappa}^{2}\right)^{-D / 2} d t = c_{\kappa} \frac{\sqrt{\pi} \Gamma\left(\frac{D-1}{2}\right)}{\Gamma\left(\frac{D}{2}\right)} = c_{\kappa} B\left(\frac{1}{2}, \frac{D-1}{2}\right)$, so the prefactor is exactly its reciprocal. Also, the expression depends only on $t^{2}$, so $f_{T}(-t)=f_{T}(t)$.  This is recognized as a Student's $t$-type density with $\nu = D - 1 = n_{1} + n_{2} - 2$ and scale $c_{\kappa}$.

We next obtain the derivative of the harmonic solution $u(r, \varphi)$ with respect to $t$.  First, recall the zonal Gegenbauer expansion on the unit ball:
\begin{eqnarray*}
u(r, \varphi)=\sum_{N=0}^{\infty} A_{N}(t) r^{N} C_{N}^{(\alpha)}(\cos \varphi), \qquad \alpha=\frac{D-2}{2}.
\end{eqnarray*}
Only the boundary aperture $\varphi_{0}(\tau)=\operatorname{arccot}(\tau)$ (with $\tau=t / c_{\kappa}$) makes $A_{N}$ depend on $t$.  From the definition, we have 
\begin{eqnarray*}
A_{N}(t) = \frac{1}{h_{N}^{(\alpha)}} \int_{\varphi_{0}(\tau)}^{\pi} C_{N}^{(\alpha)}(\cos \phi)(\sin \phi)^{2 \alpha} d \phi,
\end{eqnarray*}
Using the Leibniz rule for differentiation under the integral sign when the limits depend on the parameter (only the lower limit depends on $t$), we get
\begin{eqnarray*}
\frac{d A_{N}}{d t} &=& -\frac{1}{h_{N}^{(\alpha)}} C_{N}^{(\alpha)}\left(\cos \varphi_{0}\right)\left(\sin \varphi_{0}\right)^{2 \alpha} \frac{d \varphi_{0}}{d t} \\
&=& \frac{1}{c_{\kappa} h_{N}^{(\alpha)}} C_{N}^{(\alpha)}\left(\frac{\tau}{\sqrt{1+\tau^{2}}}\right) \left(\frac{1}{\sqrt{1+\tau^{2}}}\right)^{2 \alpha} \frac{1}{1+\tau^{2}} \\
&=& \frac{1}{c_{\kappa} h_{N}^{(\alpha)}} C_{N}^{(\alpha)}\left(c_{0}\right)\left(1+\tau^{2}\right)^{-(\alpha+1)},
\end{eqnarray*}
where $d \varphi_{0}/d \tau = -1/(1+\tau^{2})$, $d \varphi_{0}/d t = -1/\left[c_{\kappa}\left(1+\tau^{2}\right)\right]$, $\cos \varphi_{0} = \tau/\sqrt{1+\tau^{2}}$, $\sin \varphi_{0} = 1/\sqrt{1+\tau^{2}}$.  Hence, the matching series for the $t$-derivative of the harmonic solution is
\begin{eqnarray*}
\partial_{t} u(r, \varphi) = \sum_{N=0}^{\infty} \left[ \frac{1}{c_{\kappa} h_{N}^{(\alpha)}} C_{N}^{(\alpha)}\left(c_{0}\right)\left(1+\tau^{2}\right)^{-(\alpha+1)} \right] r^{N} C_{N}^{(\alpha)}(\cos \varphi),
\end{eqnarray*}
where $c_{0} = \tau / \sqrt{1+\tau^{2}}$, $\tau=t / c_{\kappa}$, $\alpha=(D-2) / 2$, and
\begin{eqnarray*}
h_{N}^{(\alpha)} = \int_{0}^{\pi}\left(C_{N}^{(\alpha)}(\cos \phi)\right)^{2}(\sin \phi)^{2 \alpha} d \phi = \frac{\pi 2^{1-2 \alpha} \Gamma(N+2 \alpha)}{N !(N+\alpha) \Gamma(\alpha)^{2}}.
\end{eqnarray*}
At $r=0$, only $N=0$ survives (since $C_{0}^{(\alpha)} \equiv 1$):
\begin{eqnarray}
\partial_{t} u(0, \cdot) = \frac{1}{c_{\kappa} h_{0}^{(\alpha)}}\left(1+\tau^{2}\right)^{-(\alpha+1)} = \frac{1}{c_{\kappa} B\left(\frac{1}{2}, a\right)}\left(1+\tau^{2}\right)^{-\frac{D}{2}} = f_{T}(t), \label{eq:pdf_final_1}
\end{eqnarray}
exactly matching the pdf in \eqref{eq:pdf_final}.  We have, therefore, expressed everything in terms of Beta-Gamma constants and standard Gegenbauer polynomials, and the apex limit of the series recovers the pdf.

\subsection{Extraction of large-$\lvert t \rvert$ asymptotics from the incomplete Beta representation}\label{subsec:tail_behavior}
Recall from \eqref{eq:cdf_spherical_cap} that the cdf in incomplete beta form is given by 
\begin{eqnarray*}
F_T(t)=\frac{1}{2}\left[1+\operatorname{sgn}(t) I_{\xi(t)}\left(\frac{1}{2}, a\right)\right],\qquad \xi(t):=\frac{\tau^{2}}{1+\tau^{2}}=\frac{t^{2}}{t^{2}+c_{\kappa}^{2}}.
\end{eqnarray*}
For $t>0$, $F_T(t)=\frac{1}{2}\left[1+I_{\xi(t)}\left(\frac{1}{2}, a\right)\right]$; for $t < 0$ we can use symmetry.  We will use the standard edge expansion, for $x \uparrow 1$ and $q>0$,
\begin{eqnarray}
1-I_x(p, q)=\frac{(1-x)^{q}}{q B(p, q)}\left[1-\frac{p-1}{q+1}(1-x)+O\left((1-x)^{2}\right)\right].   \label{eq:edge_expansion}
\end{eqnarray}
To characterize the right-tail of the cdf $1-F_T(t)$ as $t \to +\infty$, for $t > 0$, we have 
\begin{eqnarray*}
1-F_T(t)=\frac{1}{2}\left[1-I_{\xi(t)}\left(\frac{1}{2}, a\right)\right].
\end{eqnarray*}
As $t \to +\infty$, $\xi(t)=\frac{\tau^{2}}{1+\tau^{2}} \uparrow 1$, with
\begin{eqnarray*}
1-\xi(t)=\frac{1}{1+\tau^{2}}=\frac{c_{\kappa}^{2}}{t^{2}+c_{\kappa}^{2}}=\frac{c_{\kappa}^{2}}{t^{2}}\left(1+O\left(t^{-2}\right)\right).
\end{eqnarray*}
With $p=\frac{1}{2}$, $q=a$, and $x=\xi(t)$, we have 
\begin{eqnarray}
1-F_T(t) &=& \frac{1}{2} \cdot \frac{(1-\xi(t))^{a}}{a B\left(\frac{1}{2}, a\right)}\left[1-\frac{\frac{1}{2}-1}{a+1}(1-\xi(t))+O\left((1-\xi)^{2}\right)\right] \nonumber \\
&=& \frac{1}{2} \cdot \frac{(1-\xi(t))^{a}}{a B\left(\frac{1}{2}, a\right)}\left[1+\frac{1}{2(a+1)}(1-\xi(t))+O\left((1-\xi)^{2}\right)\right] \nonumber \\
&=& \frac{c_{\kappa}^{D-1}}{(D-1) B\left(\frac{1}{2}, \frac{D-1}{2}\right)} t^{-(D-1)}\left[1+\frac{c_{\kappa}^{2}}{D+1} t^{-2}+O\left(t^{-4}\right)\right], \quad t \to+\infty.  \label{eq:cdf_plus}
\end{eqnarray}
where $1-\xi(t) = c_{\kappa}^{2} / t^{2}+O\left(t^{-4}\right)$ and $2 a=D-1$.  Thus the right tail index is $D-1$ with the explicit leading constant $\frac{c_{\kappa}^{D-1}}{(D-1) B(1 / 2,(D-1) / 2)}$, and we have also obtained the next correction $O\left(t^{-(D+1)}\right)$.  By symmetry, for the left tail $t \to -\infty$, we have 
\begin{eqnarray}
F_T(t) = \frac{c_{\kappa}^{D-1}}{(D-1) B\left(\frac{1}{2}, \frac{D-1}{2}\right)} |t|^{-(D-1)}\left[1+\frac{c_{\kappa}^{2}}{D+1} t^{-2}+O\left(t^{-4}\right)\right].  \label{eq:cdf_minus}
\end{eqnarray}
Next, we derive the tail behavior of the pdf.  Differentiate the leading tail \eqref{eq:cdf_plus}
\begin{eqnarray*}
\frac{d}{d t}\left(t^{-(D-1)}\right)=-(D-1) t^{-D}, \qquad \frac{d}{d t}\left(t^{-(D+1)}\right)=-(D+1) t^{-(D+2)}.
\end{eqnarray*}
Hence, we get 
\begin{eqnarray}
f_T(t) &=& \frac{c_{\kappa}^{D-1}}{B\left(\frac{1}{2}, \frac{D-1}{2}\right)} t^{-D}\left[1-\frac{(D-1) c_{\kappa}^{2}}{D+1} t^{-2}+O\left(t^{-4}\right)\right], \quad t \to+\infty.  \label{eq:pdf_from_cdf}
\end{eqnarray}
This matches the expression in \eqref{eq:pdf_final}
\begin{eqnarray*}
f_T(t) = \frac{1}{c_{\kappa} B\left(\frac{1}{2}, \frac{D-1}{2}\right)}\left(1+\frac{t^{2}}{c_{\kappa}^{2}}\right)^{-\frac{D}{2}}
\end{eqnarray*}
expanded as
\begin{eqnarray*}
\left(1+\frac{t^{2}}{c_{\kappa}^{2}}\right)^{-D / 2}=\left(\frac{t^{2}}{c_{\kappa}^{2}}\right)^{-D / 2}\left(1+\frac{c_{\kappa}^{2}}{t^{2}}\right)^{-D / 2}=\frac{c_{\kappa}^{D}}{t^{D}}\left[1-\frac{D}{2} \frac{c_{\kappa}^{2}}{t^{2}}+O\left(t^{-4}\right)\right],
\end{eqnarray*}
and multiplying by $1 /\left(c_{\kappa} B\right)$ gives exactly the same leading constant and the first correction (since $\frac{D}{2}=\frac{D-1}{2}+\frac{1}{2}$ and the two derivations agree up to algebraic identity).

Lastly, we derive the tail asymptotics for $\partial_t u(r, \varphi)$ (the field derivative) as $t \to \infty$.  From the series expansion, we write 
\begin{eqnarray*}
\partial_t u(r, \varphi)=\sum_{N=0}^{\infty} \frac{1}{c_{\kappa} h_N^{(\alpha)}} C_N^{(\alpha)}\left(c_0\right)\left(1+\tau^{2}\right)^{-(\alpha+1)} r^{N} C_N^{(\alpha)}(\cos \varphi), \\
\tau=\frac{t}{c_{\kappa}}, \quad c_0=\frac{\tau}{\sqrt{1+\tau^{2}}}, \quad h_N^{(\alpha)}=\frac{\pi 2^{1-2 \alpha} \Gamma(N+2 \alpha)}{N!(N+\alpha) \Gamma(\alpha)^{2}}.
\end{eqnarray*}
As $t \to \infty$, $c_0 \to 1$ and $\left(1+\tau^{2}\right)^{-(\alpha+1)}=\tau^{-2(\alpha+1)}\left(1+O\left(\tau^{-2}\right)\right)$. Use the endpoint value
\begin{eqnarray*}
C_N^{(\alpha)}(1)=\frac{\Gamma(N+2 \alpha)}{\Gamma(2 \alpha) N!}
\end{eqnarray*}
and the reciprocal of $h_N^{(\alpha)}$. This yields
\begin{eqnarray*}
\frac{C_N^{(\alpha)}(1)}{h_N^{(\alpha)}}=\frac{\Gamma(N+2 \alpha)}{\Gamma(2 \alpha) N!} \times \frac{N!(N+\alpha) \Gamma(\alpha)^{2}}{\pi 2^{1-2 \alpha} \Gamma(N+2 \alpha)} = \underbrace{\frac{\Gamma(\alpha)^{2}}{\pi 2^{1-2 \alpha} \Gamma(2 \alpha)}}_{\eqqcolon K_\alpha} (N+\alpha),
\end{eqnarray*}
{\ie}, independent of $N$ up to the factor $(N+\alpha)$. Therefore,
\begin{eqnarray*}
\partial_t u(r, \varphi) = \frac{\left(1+\tau^{2}\right)^{-(\alpha+1)}}{c_{\kappa}} K_\alpha \sum_{N=0}^{\infty}(N+\alpha) r^{N} C_N^{(\alpha)}(\cos \varphi)+O\left(\left(1+\tau^{2}\right)^{-(\alpha+2)}\right).
\end{eqnarray*}
The Poisson kernel on $S^{D-1}$ (zonal form) gives the identity
\begin{eqnarray*}
\sum_{N=0}^{\infty}(N+\alpha) r^{N} C_N^{(\alpha)}(\cos \varphi) &=& \alpha \frac{1-r^{2}}{\left(1-2 r \cos \varphi+r^{2}\right)^{\alpha+1}}, \quad 0 \le r < 1. 
\end{eqnarray*}
Hence, we get
\begin{eqnarray*}
\partial_t u(r, \varphi)=\frac{\left(1+\tau^{2}\right)^{-(\alpha+1)}}{c_{\kappa}} K_\alpha \alpha \frac{1-r^{2}}{\left(1-2 r \cos \varphi+r^{2}\right)^{\alpha+1}}+O\left(\left(1+\tau^{2}\right)^{-(\alpha+2)}\right).
\end{eqnarray*}
Using the duplication formula $\Gamma(2 \alpha)=\frac{2^{2 \alpha-1}}{\sqrt{\pi}} \Gamma(\alpha) \Gamma\left(\alpha+\frac{1}{2}\right)$, we can write
\begin{eqnarray*}
K_\alpha \alpha=\frac{\alpha \Gamma(\alpha)^{2}}{\pi 2^{1-2 \alpha} \Gamma(2 \alpha)}=\frac{\alpha \Gamma(\alpha)}{\sqrt{\pi} \Gamma\left(\alpha+\frac{1}{2}\right)}=\frac{1}{B\left(\frac{1}{2}, a\right)}, \quad a=\alpha+\frac{1}{2}.
\end{eqnarray*}
Therefore, the leading tail of the field derivative is
\begin{eqnarray}
\partial_t u(r, \varphi)=\frac{1}{c_{\kappa} B\left(\frac{1}{2}, a\right)}\left(1+\tau^{2}\right)^{-(\alpha+1)} \frac{1-r^{2}}{\left(1-2 r \cos \varphi+r^{2}\right)^{\alpha+1}}+O\left(\left(1+\tau^{2}\right)^{-(\alpha+2)}\right). \label{eq:field_tail}
\end{eqnarray}
At the apex ($r=0$), the geometric factor equals 1 and \eqref{eq:field_tail} reduces to
\begin{eqnarray}
\partial_t u(0, \cdot) = \frac{1}{c_{\kappa} B\left(\frac{1}{2}, a\right)}\left(1+\tau^{2}\right)^{-(\alpha+1)}=\frac{1}{c_{\kappa} B\left(\frac{1}{2}, \frac{D-1}{2}\right)}\left(1+\frac{t^{2}}{c_{\kappa}^{2}}\right)^{-D / 2}=f_T(t), \label{eq:field_tail_pdf}
\end{eqnarray}
as required.  Keeping the next term in $C_N^{(\alpha)}\left(c_0\right)=C_N^{(\alpha)}(1)+O(1-\xi)$ with
\begin{eqnarray*}
1 - \xi = \frac{1}{1+\tau^{2}}=\frac{c_{\kappa}^{2}}{t^{2}}+O\left(t^{-4}\right),
\end{eqnarray*}
and the expansion $\left(1+\tau^{2}\right)^{-(\alpha+1)}=\tau^{-2(\alpha+1)}\left(1-(\alpha+1) \tau^{-2}+O\left(\tau^{-4}\right)\right)$, providing a systematic $O\left(t^{-(D+2)}\right)$ correction to \eqref{eq:field_tail}. The coefficients can be written explicitly by differentiating the closed form of the Poisson kernel with respect to the parameter $x=c_0$ at $x=1$; we omit that algebra here since it parallels the cdf correction already appearing in \eqref{eq:cdf_plus}.  These formulas give precise constants and first corrections, directly traced to the incomplete Beta edge $x \uparrow 1$ and the Poisson kernel on the sphere, yielding tail exponents $D-1$ for the cdf and $D$ for the pdf.

\subsection{Numerical results}\label{subsec:numerical_results}
\begin{figure}[H]
\centering
\begin{subfigure}{0.45\textwidth}
    \centering
    \includegraphics[width=\linewidth]{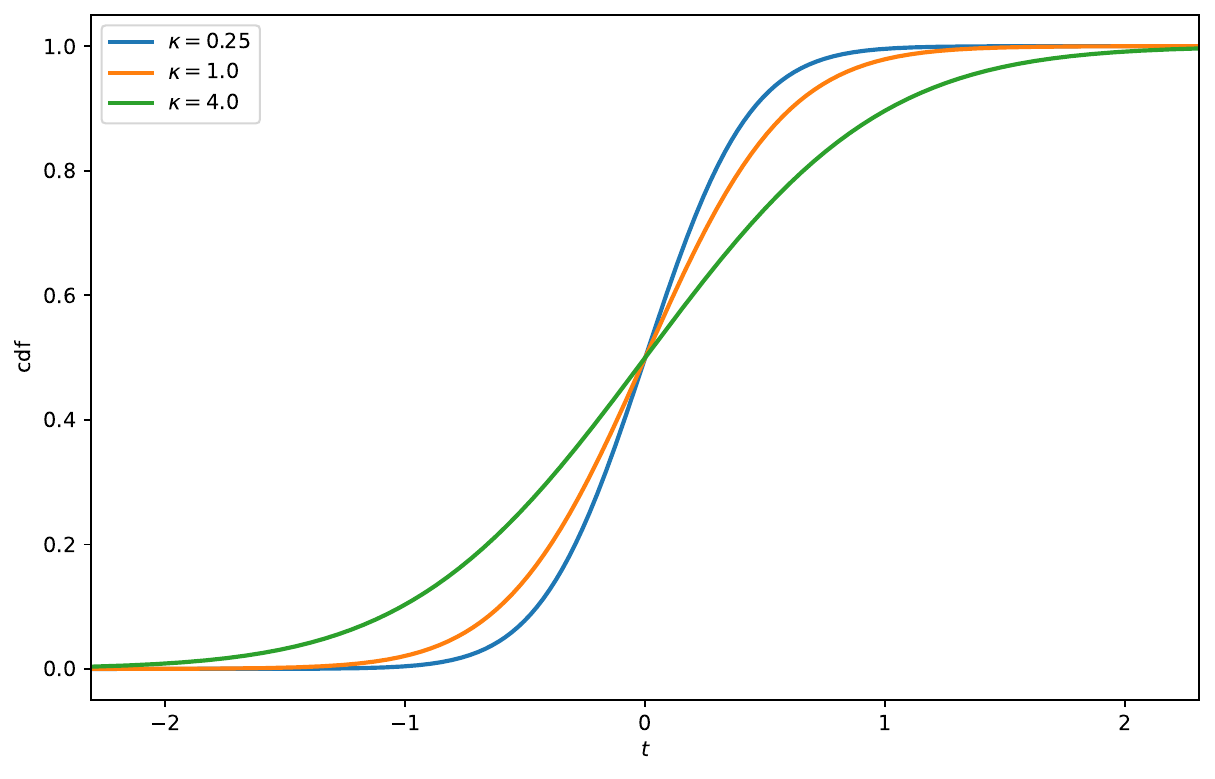}
    \caption{Effect of $\kappa$ on the cdf.}
    \label{fig:laplace_bf_cdf_variance_ratio}
\end{subfigure}
\hfill
\begin{subfigure}{0.45\textwidth}
    \centering
    \includegraphics[width=\linewidth]{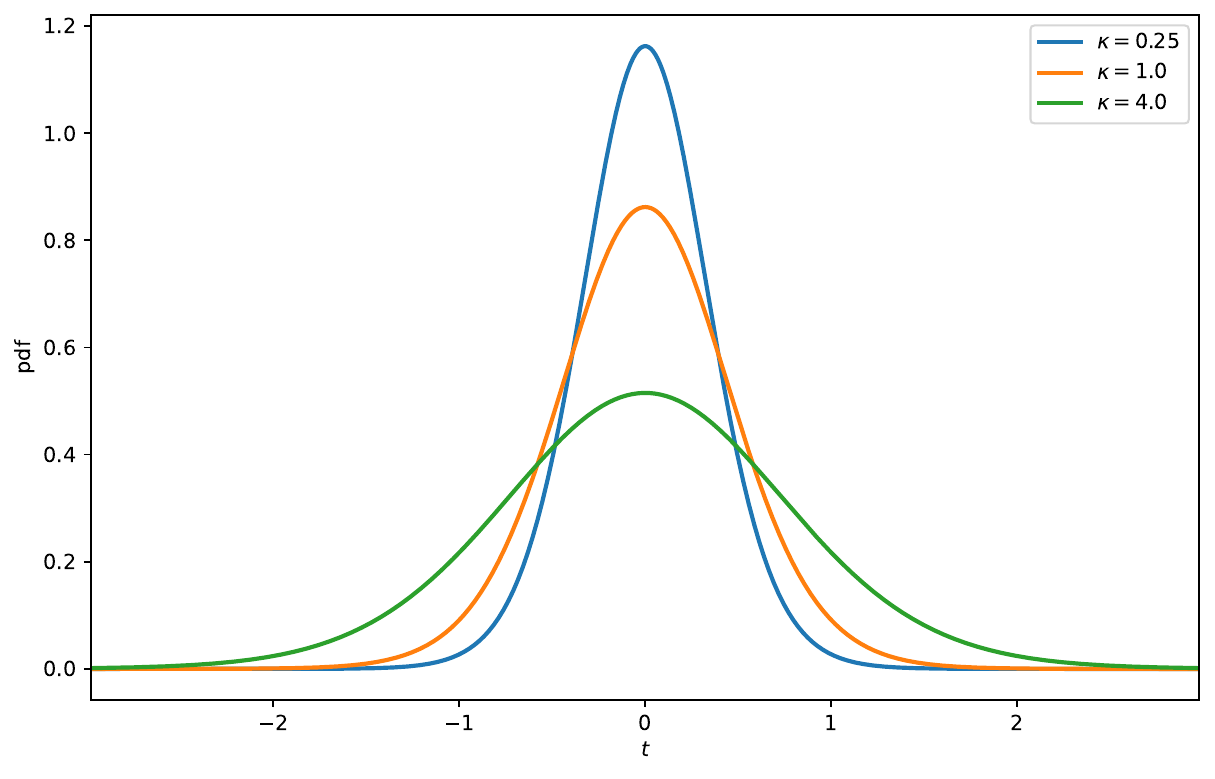}
    \caption{Effect of $\kappa$ on the pdf.}
    \label{fig:laplace_bf_pdf_variance_ratio}
\end{subfigure}
\caption{Effect of the variance ratio $\kappa$ on the distribution and density functions.}
\label{fig:kappa_cdf_pdf}
\end{figure}

\begin{figure}[H]
\centering
\begin{subfigure}{0.45\textwidth}
    \centering
    \includegraphics[width=\linewidth]{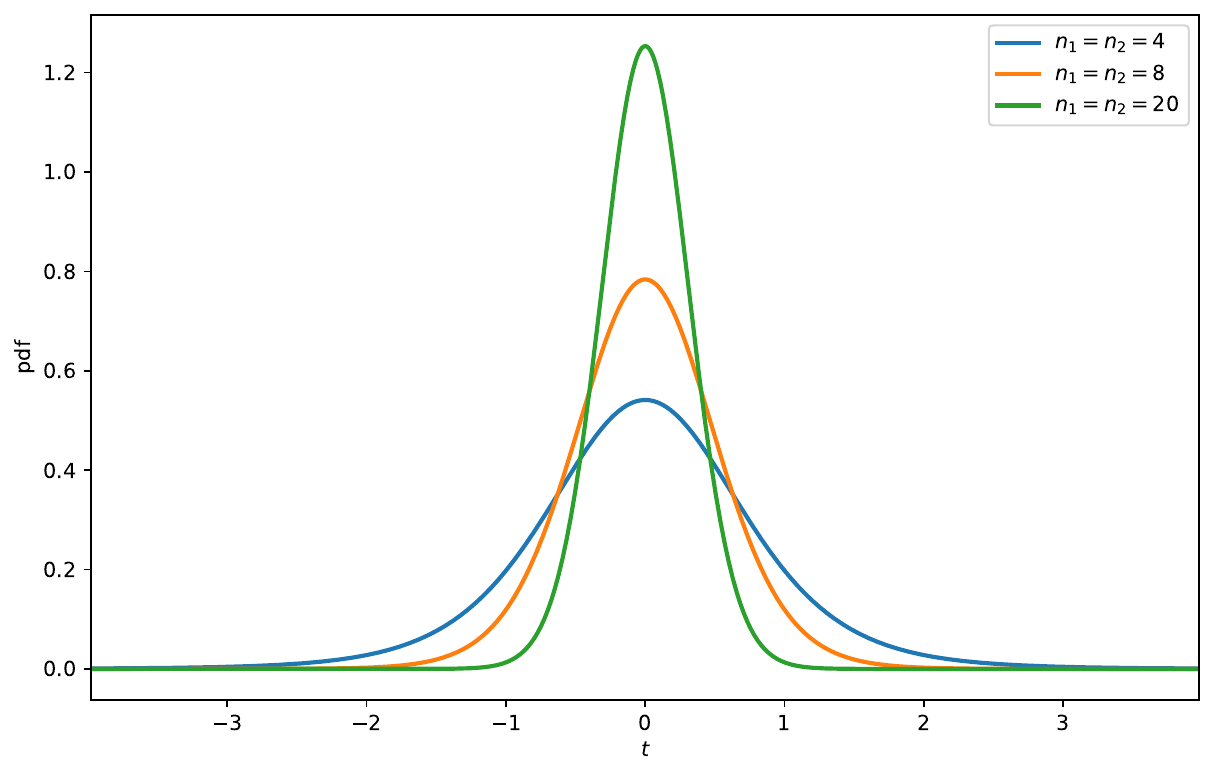}
    \caption{Effect of $(n_1, n_2)$ on the pdf.}
    \label{fig:laplace_bf_pdf_sample_size}
\end{subfigure}
\hfill
\begin{subfigure}{0.45\textwidth}
    \centering
    \includegraphics[width=\linewidth]{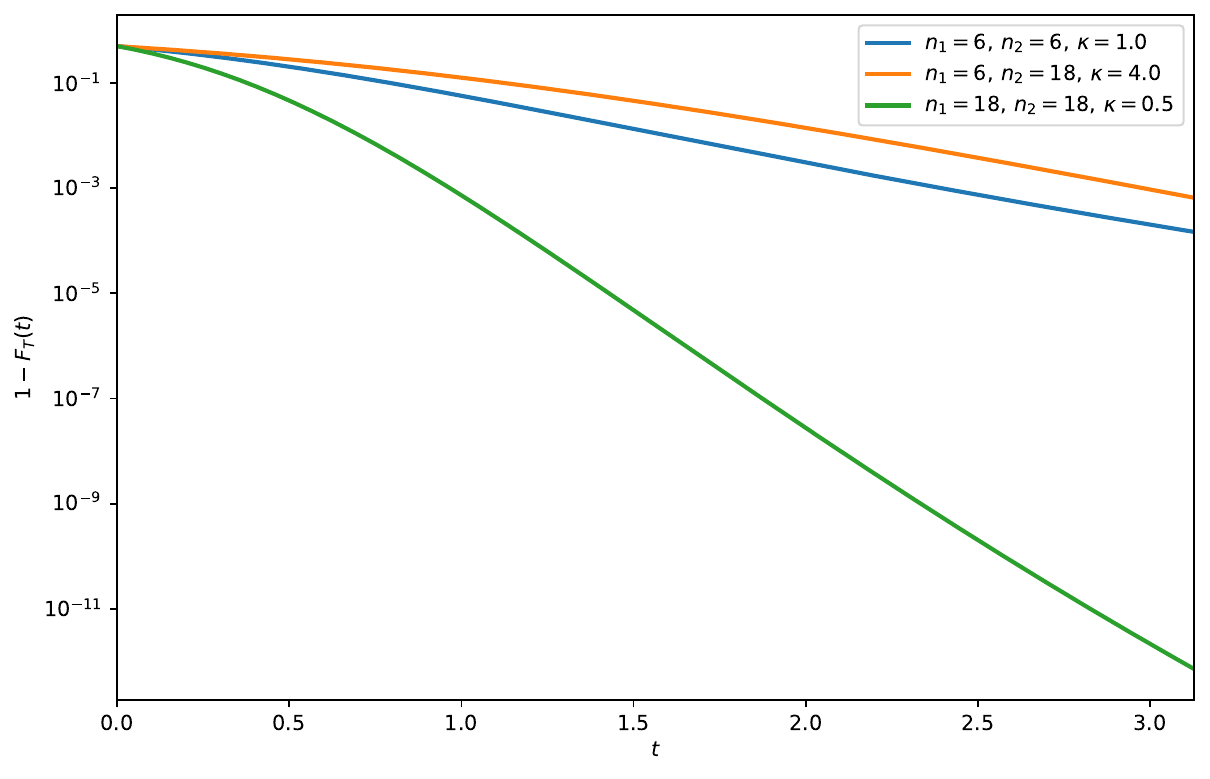}
    \caption{Effect of $(n_1, n_2, \kappa)$ on the tail behavior.}
    \label{fig:laplace_bf_tail_logscale}
\end{subfigure}
\caption{Effect of sample sizes $(n_1, n_2, \kappa)$ on the density function.}
\label{fig:sample_size_pdf}
\end{figure}

To illustrate the finite-sample behaviour of the distribution and density functions, we evaluated the closed-form expressions for the cdf in \eqref{eq:cdf_final} and the pdf in \eqref{eq:pdf_final} over a grid of threshold values and across several representative parameter configurations. The numerical study was designed to isolate the separate roles of the variance ratio $\kappa=\sigma_1^2/\sigma_2^2$ and the sample-size pair $(n_1,n_2)$.  In \figref{fig:kappa_cdf_pdf}, $(n_1,n_2)$ is held fixed while $\kappa$ varies, so that the effect of heteroscedasticity on the exact finite-sample law can be examined directly; the cdf is shown in \figref{fig:laplace_bf_cdf_variance_ratio} and the pdf is shown in \figref{fig:laplace_bf_pdf_variance_ratio}.  In \figref{fig:laplace_bf_pdf_sample_size}, $\kappa$ is fixed and the sample sizes are varied, in order to display the dependence of the distribution on $D=n_1+n_2-1$.  In \figref{fig:laplace_bf_tail_logscale}, $(n_1, n_2, \kappa)$ are varied to see the tail behavior.   Since the graphs are computed directly from the explicit beta-function and incomplete-beta-function representations, the plots reflect the exact finite-sample formulas derived here, rather than any asymptotic approximation or Monte Carlo simulations.  
The graphical behaviour is fully consistent with the analytic structure of the law. Variation in $\kappa$ acts through the scale factor $c_{\kappa}=\sqrt{\kappa/n_1+1/n_2}$, thereby altering both the dispersion of the density and the rate of transition of the distribution function as shown in \figref{fig:kappa_cdf_pdf}. For fixed variance ratio, increasing the sample sizes produces the expected concentration of the law, visible in the density as a sharper central peak and a corresponding reduction in dispersion as shown in \figref{fig:laplace_bf_pdf_sample_size}. In addition, a logarithmic tail plot of $1-F_T(t)$ shown in \figref{fig:laplace_bf_tail_logscale} provides a useful comparative view of tail decay across parameter settings and complements the explicit tail asymptotics established earlier.  The graphical results provide a concrete numerical realization of the exact theory, clarifying how heteroscedasticity and sample size jointly shape the finite-sample law while also demonstrating the computational tractability of the derived formulas.

\begin{table}[htbp]
\centering
\normalsize
\setlength{\tabcolsep}{6pt}
\renewcommand{\arraystretch}{1.2}

\begin{minipage}{0.48\textwidth}
\centering
\begin{tabular}{ccccc}
\multicolumn{5}{c}{$(n_1,n_2,\kappa)=(8,8,1.0)$} \\
\hline
$p$ & Exact & Welch & Diff. & \% over \\
\hline
0.900 & 1.3464 & 1.3450 & -0.0014 & -0.10 \\
0.950 & 1.7627 & 1.7613 & -0.0014 & -0.08 \\
0.975 & 2.1466 & 2.1448 & -0.0018 & -0.08 \\
0.990 & 2.6272 & 2.6245 & -0.0027 & -0.10 \\
\hline
\end{tabular}
\end{minipage}
\hfill
\begin{minipage}{0.48\textwidth}
\centering
\begin{tabular}{ccccc}
\multicolumn{5}{c}{$(n_1,n_2,\kappa)=(8,8,4.0)$} \\
\hline
$p$ & Exact & Welch & Diff. & \% over \\
\hline
0.900 & 1.3674 & 1.3694 & 0.0020 & 0.15 \\
0.950 & 1.8070 & 1.8073 & 0.0003 & 0.02 \\
0.975 & 2.2193 & 2.2198 & 0.0005 & 0.02 \\
0.990 & 2.7306 & 2.7494 & 0.0188 & 0.69 \\
\hline
\end{tabular}
\end{minipage}

\vspace{0.8em}
\begin{minipage}{0.48\textwidth}
\centering
\begin{tabular}{ccccc}
\multicolumn{5}{c}{$(n_1,n_2,\kappa)=(8,8,1.0)$} \\
\hline
$p$ & Exact & Welch & Diff. & \% over \\
\hline
0.900 & 1.3498 & 1.3504 & 0.0006 & 0.04 \\
0.950 & 1.7710 & 1.7714 & 0.0004 & 0.03 \\
0.975 & 2.1613 & 2.1611 & -0.0002 & -0.01 \\
0.990 & 2.6512 & 2.6515 & 0.0003 & 0.01 \\
\hline
\end{tabular}
\end{minipage}
\hfill
\begin{minipage}{0.48\textwidth}
\centering
\begin{tabular}{ccccc}
\multicolumn{5}{c}{$(n_1,n_2,\kappa)=(8,8,4.0)$} \\
\hline
$p$ & Exact & Welch & Diff. & \% over \\
\hline
0.900 & 1.3898 & 1.3903 & 0.0005 & 0.04 \\
0.950 & 1.8467 & 1.8470 & 0.0004 & 0.02 \\
0.975 & 2.2846 & 2.2850 & 0.0004 & 0.02 \\
0.990 & 2.8570 & 2.8605 & 0.0036 & 0.12 \\
\hline
\end{tabular}
\end{minipage}
\caption{Our results and Welch quantiles for selected configurations.}
\label{tab:quantiles_comparison}
\end{table}
For each selected triple $(n_1,n_2,\kappa)$, the critical values were computed at fixed nominal upper-tail probabilities using two parallel procedures.  The results are tabulated in \tabref{tab:quantiles_comparison}.  Exact critical values were obtained from the finite-sample distribution law derived in this paper.  The exact quantile is given by $q_{\text{exact}}(p) = \sqrt{\nu} c_{\kappa} t_{\nu, p}$, $\nu = n_1 + n_2 - 2$, $c_\kappa = \sqrt{\kappa/n_1 + 1/n_2}$, where $t_{\nu, p}$ is the $p^{\mathrm{th}}$ quantile of the standard Student $t$ law with $\nu$ degrees of freedom.  

For the Welch comparison, the $p^{\text{th}}$ quantile was computed by $200{,}000$ Monte Carlo replications at each fixed parameter configuration $(n_1,n_2,\kappa)$. Specifically, for each replicate we generated one pair of normal samples under the prescribed variances, computed the corresponding sample variances $S_1^2$ and $S_2^2$, substituted these values into the Welch--Satterthwaite formula to obtain a replicate-specific effective degrees of freedom $\nu_W$, and then evaluated the Student $t$ quantile $q_{\text{Welch}}(p) = t_{\nu_W, p}$, where
\begin{eqnarray*}
\nu_W = \frac{\left( \frac{\kappa}{n_1} + \frac{1}{n_2} \right)^2}{\frac{\kappa^2}{n_1^2(n_1 - 1)} + \frac{1}{n_2^2(n_2 - 1)}}.
\end{eqnarray*}
Repeating this procedure produced $200{,}000$ replicate-specific Welch critical values, and the reported Welch entry in the table was taken to be their Monte Carlo average. Thus, the Welch quantiles shown in \tabref{tab:quantiles_comparison} are not based on a single fixed degrees-of-freedom approximation, but on an average over the random finite-sample variability induced by the sample variances themselves.  \tabref{tab:quantiles_comparison} indicates that, for the configurations reported, the Welch critical values are numerically very close to the compact finite-sample quantiles derived here. This demonstrates that the cdf and pdf laws developed in this paper provides the benchmark that makes such a comparison possible, and it yields a rigorous analytic characterization of the Behrens--Fisher statistic that the Welch approximation does not provide. In particular, the PDE-based derivation identifies the precise finite-sample distribution, not merely an effective degrees-of-freedom surrogate, and it supplies additional structural information, including closed-form distributional representations and tail behaviour, beyond what can be read off from the classical approximation. 

\section{Conclusion}\label{sec:conclusion} \vspace{-0.5cm}
A Laplace-equation framework was developed for the Behrens--Fisher problem under unknown and unequal variances by reducing the event defining the statistic to a scale-invariant conic inequality and then identifying the associated spherical probability with harmonic measure. This converts the finite-sample distributional problem into a Laplace--Dirichlet boundary value problem on the unit ball and leads to explicit representations for the cumulative distribution function and density in terms of regularized incomplete beta functions, together with a Gegenbauer separation-of-variables expansion for the associated harmonic extension and its threshold derivative. The analysis also yields sharp tail expansions with explicit leading constants and higher-order corrections, thereby providing both structural insight and practically computable formulas.  Beyond the present normal-theory setting, the geometric formulation suggests several natural directions for further work, including extensions to more general elliptically contoured models, analogous boundary-value formulations for other heteroscedastic studentized statistics, and a deeper spectral analysis of the wedge problem to obtain refined eigenvalue asymptotics and rigorous truncation-error bounds for the associated series representations.  The quantile tables indicate that, in the representative settings considered, the exact finite-sample critical values are in close proximity to the Welch--Satterthwaite approximation. This suggests that the exact law derived in this paper can provide useful calibration beyond what is available from a single effective degrees-of-freedom approximation.

\section*{Disclosure of interest}
The authors report there are no competing interests to declare.

\bibliography{/Users/kgnagananda/Documents/Work/collaborations/pdx/research/references/research_pdx.bib}

\end{document}